# EMPIRICAL BAYES SELECTION OF WAVELET THRESHOLDS

By Iain M. Johnstone[1] and Bernard W. Silverman[2]

*Stanford University and University of Oxford*

This paper explores a class of empirical Bayes methods for level-dependent threshold selection in wavelet shrinkage. The prior considered for each wavelet coefficient is a mixture of an atom of probability at zero and a heavy-tailed density. The mixing weight, or sparsity parameter, for each level of the transform is chosen by marginal maximum likelihood. If estimation is carried out using the posterior median, this is a random thresholding procedure; the estimation can also be carried out using other thresholding rules with the same threshold. Details of the calculations needed for implementing the procedure are included. In practice, the estimates are quick to compute and there is software available. Simulations on the standard model functions show excellent performance, and applications to data drawn from various fields of application are used to explore the practical performance of the approach.

By using a general result on the risk of the corresponding marginal maximum likelihood approach for a single sequence, overall bounds on the risk of the method are found subject to membership of the unknown function in one of a wide range of Besov classes, covering also the case of $f$ of bounded variation. The rates obtained are optimal for any value of the parameter $p$ in $(0,\infty]$, simultaneously for a wide range of loss functions, each dominating the $L_q$ norm of the $\sigma$th derivative, with $\sigma \geq 0$ and $0 < q \leq 2$.

Attention is paid to the distinction between sampling the unknown function within white noise and sampling at discrete points, and between placing constraints on the function itself and on the discrete wavelet transform of its sequence of values at the observation points. Results for all relevant combinations of these scenarios are obtained.

Received September 2003; revised September 2004.

[1]Supported in part by NSF Grant DMS-00-72661 and NIH Grants California 72028 and R01 EB001988-08.

[2]Supported in part by by the Engineering and Physical Sciences Research Council.

*AMS 2000 subject classifications.* Primary 62C12, 62G08; secondary 62G20, 62H35, 65T60.

*Key words and phrases.* Adaptivity, Bayesian inference, nonparametric regression, smoothing, sparsity.

> In some cases a key feature of the theory is a particular boundary-corrected wavelet basis, details of which are discussed.
>
> Overall, the approach described seems so far unique in combining the properties of fast computation, good theoretical properties and good performance in simulations and in practice. A key feature appears to be that the estimate of sparsity adapts to three different zones of estimation, first where the signal is not sparse enough for thresholding to be of benefit, second where an appropriately chosen threshold results in substantially improved estimation, and third where the signal is so sparse that the zero estimate gives the optimum accuracy rate.

## 1. Introduction.

1.1. *Background.* Consider the nonparametric regression problem where we have observations at $2^J$ regularly spaced points $t_i$ of some unknown function $f$ subject to noise

$$X_i = f(t_i) + \varepsilon_i, \tag{1}$$

where the $\varepsilon_i$ are independent $N(0, \sigma_E^2)$ random variables. The standard wavelet-based approaches to the estimation of $f$ proceed by taking the discrete wavelet transform of the data $X_i$, processing the resulting coefficients to remove noise, usually by some form of thresholding, and then transforming back to obtain the estimate.

The underlying notion behind wavelet methods is that the unknown function has an economical wavelet expression, in that $f$ is, or is well approximated by, a function with a relatively small proportion of nonzero wavelet coefficients. The quality of estimation is quite sensitive to the choice of threshold, with the best choice being dependent on the problem setting. In general terms, "sparse" signals call for relatively high thresholds ($3\sigma_E$, $4\sigma_E$ or even higher), while "dense" signals might demand choices of $2\sigma_E$ or even lower. Indeed, it is typical that the wavelet coefficients of a true signal will be relatively more sparse at the fine resolution scales than at the coarser scales. It is therefore desirable to develop threshold selection methods that adapt the threshold level by level.

One would hope that such methods would estimate thresholds that stably reflect the gradation from sparse to dense signals as the scale changes from fine to coarse. It has proven elusive to construct threshold selectors that combine properties such as these with good theoretical properties. The principal motivation for the work reported in the present paper is to show that a simple empirical Bayesian approach combines computational tractability with good theoretical and practical performance. For software availability, see Section 1.8.

While the present paper is concerned with the nonparametric regression model (1) and wavelet transforms, the same levelwise empirical Bayes approach is, in principle, directly applicable to other direct and indirect transform shrinkage settings with multiscale block structure, as briefly discussed in [28].

1.2. *Bayesian approaches to wavelet regression.* Within a Bayesian context, the notion of sparsity is naturally modeled by a suitable prior distribution for the wavelet coefficients of $f$. Write $d_{jk}$ for the elements of the discrete wavelet transform (DWT) of the vector of values $f(t_i)$ and $d^*_{jk}$ for the DWT of the observed data $X_i$. Let $N = 2^J$ and $\theta_{jk} = N^{-1/2} d_{jk}$.

Clyde, Parmigiani and Vidakovic [13], Abramovich, Sapatinas and Silverman [4] and Silverman [46] have considered a particular mixture prior for this problem. Under this prior, the $d_{jk}$ are independently distributed with

$$(2) \qquad d_{jk} \sim (1 - \pi_j)\delta(0) + \pi_j N(0, \tau_j^2),$$

a mixture of an atom of probability at zero and a normal distribution with variance $\tau_j^2$. The parameters of the distribution (2) depend on the level $j$ of the coefficient in the transform. A related prior was considered by Chipman, Kolaczyk and McCulloch [11]; for a survey of work in this area, see [48]. See also [12, 38, 44, 50, 51] for a range of approaches to the modeling of the wavelet coefficients underlying a function or image. [31] is an early version introducing the approach of the present paper.

The most popular summary of the posterior distribution under the model (2) has been the posterior mean, but Abramovich, Sapatinas and Silverman [4] investigated the use of the posterior median of $d_{jk}$ as a summary of the posterior distribution. This is a true thresholding rule, in that for $|d^*_{jk}|$ less than some threshold, the point estimate of $d_{jk}$ will be exactly zero. In the wavelet context, the coefficient-wise posterior median corresponds to a point estimate of the posterior distribution under a family of loss function equivalent to $L^1$ norms on the function and its derivatives. Such $L^1$ losses are in any case more natural if one wishes to allow for the possibility of inhomogeneous functions, one of the aims of the wavelet approach.

1.3. *Choosing the parameters in the prior.* How should the parameters in the prior be chosen? In much of the existing literature, the parameters are either chosen directly by reference to prior information about $f$, or by a combination of prior information and data-based criteria. Though some of these, for example, the BayesThresh approach of Abramovich, Sapatinas and Silverman [4], give good results, they clearly invite the possibility of a more systematic approach to the choice of the hyperparameters. In the present paper we take an empirical Bayes (or marginal maximum likelihood) approach, which yields a completely data-based method of choosing the prior parameters. Within the Bayesian formulation set out above, wavelet regression at a single resolution level $j$ is a special case of a single sequence Bayesian model selection problem considered, among others, by George and Foster [24, 25]. This problem is considered in detail by Johnstone and Silverman [33]; we review the basic method presented there and also give some additional implementational details.

Suppose that $Z = (Z_1, \ldots, Z_n)$ are observations satisfying

$$Z_i = \mu_i + \varepsilon_i, \tag{3}$$

where the $\varepsilon_i$ are independent $N(0,1)$ random variables. It is supposed that the unknown coefficients $\mu_i$ are mostly zero, but some of them may be nonzero, and, with this in mind, it is of interest to estimate the $\mu_i$ on the basis of the observed data. In the model selection context, the nonzero $\mu_i$ correspond to parameters that actually enter the model. The connection with wavelet regression is natural: the $Z_i$ might be the sample wavelet coefficients (suitably renormalized) at a particular level, and these are noisy observations of a sequence of population wavelet coefficients which are mostly zero.

The parameters $\mu_i$ are modeled as having independent prior distributions each given by the mixture

$$f_{\text{prior}}(\mu) = (1-w)\delta_0(\mu) + w\gamma(\mu). \tag{4}$$

The nonzero part of the prior, $\gamma$, is assumed to be a fixed unimodal symmetric density. In most of the previous wavelet work cited above, the density $\gamma$ is a normal density, but we use a heavier-tailed prior, replacing the $N(0, \tau_j^2)$ part of the mixture (2) by, for example, a double exponential distribution with a scale parameter that may depend on the level of the coefficient in the transform. Another possible prior, with still heavier tails, is introduced in Section 2. Apart from the theoretical advantages of such an approach, Wainwright, Simoncelli and Willsky [50] argue that the marginal distribution of the wavelet coefficients of images arising in practice typically has tails heavier than Gaussian. In the Bayesian setup, the noise $(\varepsilon_i)$ is independent of the wavelet coefficients.

Let $g = \gamma \star \varphi$, where $\star$ denotes convolution. To avoid confusion with the scaling function of the wavelet family, we use $\varphi$ to denote the standard normal density. The marginal density of the observations $Z_i$ will then be

$$(1-w)\varphi(z) + wg(z).$$

We define the marginal maximum likelihood estimator $\hat{w}$ of $w$ to be the maximizer of the marginal log-likelihood

$$\ell(w) = \sum_{i=1}^{n} \log\{(1-w)\varphi(Z_i) + wg(Z_i)\}, \tag{5}$$

subject to the constraint on $w$ that the threshold satisfies $t(w) \leq \sqrt{2\log n}$. This upper limit on the threshold is the *universal threshold*, which has the property that it is asymptotically the largest absolute value for observations obtained from a zero signal, and can therefore be considered to be the appropriate limiting threshold as $w \to 0$.

Our basic approach is then to plug the value $\hat{w}$ back into the prior and then estimate the parameters $\mu_i$ by a Bayesian procedure using this value of $w$. Suppose $\mu$ has prior (4) and that we observe $Z \sim N(\mu, 1)$. Let $\hat{\mu}(z; w)$ be the median of the posterior distribution of $\mu$ given $Z = z$ and $\tilde{\mu}(z; w)$ its mean. If the posterior median is used, then $\mu_i$ will be estimated by $\hat{\mu}_i = \hat{\mu}(Z_i, \hat{w})$, while the corresponding estimate using the posterior mean is $\tilde{\mu}_i = \tilde{\mu}(Z_i; \hat{w})$.

For fixed $w < 1$, the function $\hat{\mu}(z; w)$ will be a monotonic function of $z$ with the thresholding property, in that there exists $t(w) > 0$ such that $\hat{\mu}(z; w) = 0$ if and only if $|z| \leq t(w)$. The estimated weight $\hat{w}$ thus yields an estimated threshold $t(\hat{w}) = \hat{t}$, say. A simple extension of the method is to retain the threshold $\hat{t}$ but to use a more general thresholding rule, for example, hard or soft thresholding. The main emphasis of this paper is on the choice of the threshold, rather than on the choice between different thresholding rules.

The posterior mean rule $\tilde{\mu}(z; w)$ fails to have the thresholding property, and, hence, produces estimates in which, essentially, all the coefficients are nonzero. Nevertheless, it has shrinkage properties that allow it to give good results. We shall see that, both in theory and in simulation studies, the performance of the posterior mean is good, but not quite as good as the posterior median.

The same approach can be used to estimate other parameters of the prior. In particular, if a scale parameter $a$ is incorporated by considering a prior density $(1-w)\delta_0(\mu) + wa\gamma(a\mu)$, define $g_a$ to be the convolution of $a\gamma(a\cdot)$ with the normal density. Then both $a$ and $w$ can be estimated by finding the maximum over both parameters of

$$\ell(w, a) = \sum_{i=1}^{n} \log\{(1-w)\varphi(Z_i) + wg_a(Z_i)\}.$$

In the case where there is no scale parameter to be estimated, $\ell'(w)$ is a monotonic function of $w$, so its root is very easily found numerically, provided the function $g$ is tractable. If one is maximizing over both $w$ and $a$, then a package numerical maximization routine that uses gradients has been found to be an acceptably efficient way of maximizing $\ell(w, a)$.

Details of relevant calculations for some particular priors are given in Section 2.2. All these calculations are implemented in the authors' package, Johnstone and Silverman [34], and the documentation of that package gives further details beyond those given in this paper.

1.4. *Marginal maximum likelihood in the wavelet context.* In the wavelet context, the MML approach is applied to each level of the wavelet transform separately, to yield values of $w$ and, if appropriate, $a$ that depend on the level of the transform. Let $\sigma_j^2$ be the standard deviation of the noise at level $j$. Assuming that the original noise is independent, the variance $\sigma_j^2$ will be the same for all $j$ and can, as is conventional, be estimated from the median of the absolute values of the coefficients at the highest level. More generally, for example, in the case of stationary correlated noise, it may be appropriate to estimate $\sigma_j$ separately for each level, at least at the higher levels of the transform. In this paper we have not considered the effect of sampling variability in the estimation of the noise variance, but that would be an interesting topic for future research.

At level $j$, define the sequence $Z_k = d_{jk}^*/\sigma_j$, and apply the single sequence MML approach to this sequence to obtain $\hat{w}_j$ and, if appropriate, estimates of any other parameters of the prior. The estimated wavelet coefficients of the discrete wavelet transform of the sequence $f(t_i)$ are then given by

$$\hat{d}_{jk} = \sigma_j \hat{\mu}(d_{jk}^*/\sigma_j; \hat{w}_j). \tag{6}$$

Assuming, without loss of generality, that the function $f$ is defined on the interval $[0,1]$ and the values $t_i = i/N$, crude estimates of the wavelet coefficients of the function $f$ are then $\hat{\theta}_{jk} = N^{-1/2} \hat{d}_{jk}$, neglecting boundary issues for the moment.

*Straightforward generalizations.* Natural generalizations of (6) include the inclusion of estimates of other parameters in the prior, as well as the use of the posterior mean instead of the posterior median, or the use of a more general thresholding rule than the posterior median, but still using the posterior median threshold $t(\hat{w})$. In addition, we consider two further generalizations.

*Modified thresholds for the estimation of derivatives.* When wavelet methods are used to estimate derivatives, it was shown by Abramovich and Silverman [5] that the appropriate universal threshold is not $\sqrt{2 \log n}$, but is a multiple of this quantity. We develop theory below using, for the estimation of derivatives, a modified threshold $t_A(w)$ given, for some appropriately chosen $A > 0$, by

$$t_A(w) = \begin{cases} t(w), & \text{if } t(w)^2 \leq 2\log n - 5 \log\log n, \\ \sqrt{2(1+A)\log n}, & \text{otherwise.} \end{cases} \tag{7}$$

*The translation-invariant wavelet transform.* It is by now well recognized that the translation-invariant wavelet transform [15], in general, gives much better results than the conventional transform applied with a fixed origin.

At each level $j$, the translation-invariant transform gives a sequence of $2^J$ values that are not actually independent. Each subsequence obtained by regular selection at intervals $2^{J-j}$ *will* be independent, and corresponds to the coefficients at level $j$ of the standard wavelet expansion with a particular choice of origin.

One way of proceeding would be to apply the empirical Bayes method entirely separately for each of these subsequences to obtain estimates of the relevant coefficients in the translation-invariant wavelet transform. It is simpler and more natural, however, to use the same estimates of the mixture hyperparameters for every position of the time origin, thereby borrowing strength in the estimation of the hyperparameters between the different positions of the origin. To obtain a single estimate at each level, we maximize the average, over choice of origin, of the marginal log-likelihood functions. This average is $2^{-(J-j)}$ times the "as-if-independent" log-likelihood function obtained by simply summing the log-likelihoods for each of the $2^J$ coefficients at level $j$ in the translation-invariant transform.

The estimates of the mixture parameters are then used to give individual posterior medians of each of the coefficients of the translation-invariant transform, and the estimated function is found by the average basis approach. Apart from the combination of log-likelihoods involved in the estimation of the hyperparameters, the translation-invariant method gives the result of applying the standard method at every possible choice of time origin, and then averaging over the position of the time origin.

Using an as-if-independent likelihood at each level to choose the hyperparameters is reminiscent of the *independence estimating equation* approach of Liang and Zeger [35] to parameter fitting in the marginal distribution of a sequence of identically distributed but nonindependent observations. Their paper was concerned with observations with generalized linear model dependence on the parameters and covariates. Because, for different choices of origin, the prior distributions on the coefficients are not, in general, generated from a single underlying prior model for the curve, our translation-invariant procedure involves a separate modeling of the prior information at each origin position, modulo $2^{J-j}$ for the coefficients at level $j$. Independence estimating equations, as we have used them, are a method of combining the separate problems of choosing the prior into a single problem at each level.

1.5. *Theoretical approach and results.* By now a classic way to study the adaptivity of wavelet smoothing methods is through the study of the worst behavior of a method when the wavelet coefficients of the function $f$ are constrained to lie in a particular Besov sequence space, corresponding to Besov function space membership of the function itself. Besov spaces are a flexible family that, depending on their parameters, can allow for varying degrees of inhomogeneity, as well as smoothness in the functions that they contain. Some relations between Besov spaces and spaces defined by $L_p$ norms on function and their derivatives are reviewed in Section 5.6. We shall show that the empirical Bayes method with a suitable function $\gamma$ automatically achieves the best possible minimax rate over a wide range of Besov spaces, including those with very low values of the parameter $p$ that allows for inhomogeneity in the unknown function $f$.

A particular case of the theory we develop is as follows; fuller details of the assumptions will be given later in the paper. Suppose that we have observations $X_i = f(t_i) + \varepsilon_i$ of a function $f$ at $N$ regularly spaced points $t_i$, with $\varepsilon_i$ independent $N(0, \sigma_E^2)$ random variables. Let $d_{jk} = N^{1/2} \theta_{jk}$ be the coefficients of an orthogonal discrete wavelet transform of the sequence $f(t_i)$, and let $d_j$ denote the vector with elements $d_{jk}$ as $k$ varies.

Assume that the coarsest level to which the wavelet transform is carried out is a fixed level $L \geq 0$. Denote by $d_{L-1}$ the vector of scaling coefficient(s) at this level. If periodic boundary conditions are being used and $N$ is a power of 2, the vector $d_j$ is of length $2^j$ if $j \geq L$ and $2^L$ if $j = L - 1$, and $N = 2^J$, where $J - 1$ is the finest level at which the sample coefficients are defined.

To allow for discrete wavelet transforms based on other boundary conditions and with values of $N$ that are other suitable multiples of powers of 2, we shall make the milder assumptions that $d_j$ is defined for $L - 1 \leq j < J$, with $L$ fixed and $J \to \infty$ as $N \to \infty$, that the sum of the lengths of the $d_j$ is equal to $N$, and that the length of each $d_j$ for $j \geq L$ is in the interval $[2^{j-1}, 2^j]$. The length of the vector $d_{L-1}$ of scaling coefficients is assumed to lie in $[2^{L-1}, 2^L]$, so that $2^{J-1} \leq N \leq 2^J$.

Estimate the coefficients $d_{jk}$ for $j \geq L$ by the estimate in (6), applying an empirical Bayes approach level by level, based on a mixture prior with a heavy-tailed nonzero component $\gamma$. The estimator can be either the posterior median or some other thresholding rule using the same threshold (and obeying a bounded shrinkage condition set out later). The scaling coefficients $d_{L-1}$ are estimated by their observed values $d^*_{L-1}$. To obtain the estimates $\hat{f}(t_i)$ of the function values $f(t_i)$, apply the inverse discrete wavelet transform to the estimated array $\hat{d}_{jk}$.

For $0 < p \leq \infty$ and $\alpha > \frac{1}{p} - \frac{1}{2}$, let $a = \alpha - \frac{1}{p} + \frac{1}{2}$. Define the Besov sequence space $b^\alpha_{p,\infty}(C)$ to be the set of all coefficient arrays $\theta$ such that

(8) $$\sum_k |\theta_{jk}|^p < C^p 2^{-apj} \quad \text{for all } j \text{ with } L - 1 \leq j < J.$$

Our theory shows that, for some constant $c$, possibly depending on $p$ and $\alpha$ but not on $N$ or $C$,

$$\sup_{\theta \in b^\alpha_{p,\infty}(C)} N^{-1} E \sum_{i=1}^N \{\hat{f}(t_i) - f(t_i)\}^2$$

(9)
$$\leq c\{C^{2/(2\alpha+1)}N^{-2\alpha/(2\alpha+1)} + N^{-1}(\log N)^4\}.$$

For fixed $C$, the second term in the bound (9) is negligible, and the rate $O(N^{-2\alpha/(2\alpha+1)})$ of decay of the mean square error is the best that can be attained over the relevant function class. The result (9) thus shows that, apart from the $O(N^{-1}\log^4 N)$ term, our estimation method simultaneously attains the optimum rate over a wide range of function classes, thus automatically adapting to the regularity of the underlying function. Under conditions we shall discuss, the Besov sequence space norm used in (8) is equivalent to a Besov function space norm on $f$ with the same parameters.

The main theorem of the paper goes considerably beyond (9), in the following respects:

- It demonstrates the optimal rate of convergence for mean $q$-norm errors for all $0 < q \leq 2$, not just the mean square error considered in (9).
- Beyond the posterior median, any thresholding method satisfying certain mild conditions can be used, and, for $1 < q \leq 2$, the results also hold for the posterior mean.
- If an appropriate modified threshold method is used, the optimality also extends to the estimation of derivatives of $f$.

Most of the existing statistical wavelet literature concentrates explicitly or implicitly on the *white noise* model, where we assume that we have independent observations of the wavelet coefficients of the function up to some resolution level. Little attention has been paid to the errors possibly introduced by the discretization of $f$. However, Donoho and Johnstone [20] discuss a form of discretization somewhat different from simple sampling at discrete points. Another issue not considered in detail in much of the present literature is the careful treatment in a statistical context of boundary-corrected wavelet methods, such as those introduced by Cohen, Daubechies and Vial [14]. In the current paper we do consider the effects of discretization and of boundary correction, and we prove theorems for both the white noise model and for a sampled data model.

In particular, suppose that the function $f$ is observed on $[0,1]$ at a regular grid of $N = 2^J$ points, subject to independent $N(0,\sigma_E^2)$ noise. Proceeding as above, but with an appropriate preconditioning of the data near the boundaries and treatment of the boundary wavelet coefficients, construct an estimate of $f$ itself by setting $\hat{f} = \sum_k \hat{\theta}_{L-1,k}\phi_{Lk} + \sum_{L \leq j < J}\sum_k \hat{\theta}_{jk}\psi_{jk}$, where $\phi_{jk}$ and $\psi_{jk}$ are the scaling functions and wavelets at scale $j$. Let $\mathcal{F}(C)$ be the class of functions $f$ whose true wavelet coefficients fall in $b_{p,\infty}^\alpha(C)$. Under appropriate mild conditions, a special case of our theory demonstrates that

(10) $$\sup_{f\in\mathcal{F}(C)} E\int_0^1 \{\hat{f}(t) - f(t)\}^2 \leq cC^{2/(2\alpha+1)}N^{-2\alpha/(2\alpha+1)} + o(N^{-2\alpha/(2\alpha+1)}).$$

Our results go far beyond mean integrated square error and consider accuracy of estimation in Besov sequence norms on the wavelet coefficients that imply good estimation of derivatives, as well as the function itself, and allow for losses in $q$-norms for $0 < q \leq 2$.

1.6. *Alternative approaches and related bibliography.* Finding a numerically simple and stable adaptive method for threshold choice with good theoretical and practical properties has proven to be elusive. A plethora of methods for choosing thresholds has been proposed (see, e.g., [49], Chapter 6). Apart from empirical Bayes methods, we note two other methods which have been accompanied by some theoretical analysis of their properties and for which software can easily be written. In both cases we set $Z_k = X_k/\sigma_E$, so that the thresholds are expressed on a renormalized scale.

Stein's Unbiased Risk Estimate (SURE) aims to minimize the mean squared error of soft thresholding, and is another method intended to be adaptive to different levels of sparsity. The threshold $\hat{t}_{\text{SURE}}$ is chosen as the minimizer (within the range $[0, \sqrt{2 \log n}\,]$) of

$$(11) \qquad \hat{U}(t) = n + \sum_{k=1}^{n} Z_k^2 \wedge t^2 - 2 \sum_{k=1}^{n} I\{Z_k^2 \leq t^2\}.$$

This does, indeed, have some good theoretical properties [19], but the same theoretical analysis, combined with simulation and practical experience, shows that the method can be unstable [19, 33] and that it does not choose thresholds well in sparse cases.

The False Discovery Rate (FDR) method is derived from the principle of controlling the false discovery rate in simultaneous hypothesis testing [7] and has been studied in detail in the estimation setting [3]. Order the data by decreasing magnitudes: $|Z|_{(1)} \geq |Z|_{(2)} \geq \cdots \geq |Z|_{(n)}$, and compare to a *quantile boundary*: $t_k = z(q/2 \cdot k/n)$, where the false discovery rate parameter $q \in (0, \frac{1}{2}]$. Define a crossing index $\hat{k}_F = \max\{k : |Z|_{(k)} \geq t_k\}$, and use this to set the threshold $\hat{t}_F = t_{\hat{k}_F}$. Although FDR threshold selection adapts very well to sparse signals [3], it does less well on dense signals of moderate size.

Overall, we shall see that empirical Bayes thresholding has some of the good properties of both SURE and FDR thresholding and deals with the transition between sparse and dense signals in a stable manner. A detailed discussion of theoretical comparisons between the various estimators is provided in Section 5.7.

1.7. *Structure of the paper.* In Section 2 we discuss various aspects of the mixture priors used later in the paper. The priors themselves are specified, and details given of formulas needed for the Bayesian calculations in practice.

We take the opportunity to give additional practical details not included in [33]. In the next two sections the practical performance of the proposed method is investigated, by simulation in Section 3, and by applications to data sets arising in practice in Section 4.

Section 5 contains the theoretical core of the paper for estimation of coefficient arrays under Besov sequence norm constraints. First, a wide-ranging result, Theorem 1, for the white noise model is stated. We then explore aspects of the boundary wavelet construction, including ways of mapping data to scaling function coefficients at the finest level. This allows for the definition of a boundary-corrected empirical Bayes estimator for the sampled data problem on a finite interval. The result we state about this estimator, Theorem 2, shows that it essentially attains the same performance as the estimator for the white noise case. Finally, the correspondences between Besov sequence and function norms are set out, specifically addressing wavelets and functions on a bounded interval. For $0 < q \leq 2$, we relate risk measures expressed in terms of wavelet coefficients to $q$-norms of appropriate derivatives.

Section 6 contains the proofs of the main theorems, starting by reviewing theoretical results for the single sequence problem from [33], but cast into a form relevant for the present paper. These results are used to prove the white noise case theorem. The proof of the theorem for the sampled data case also makes use of approximation results for appropriate boundary-corrected wavelets given by Johnstone and Silverman [32]. Finally, Section 7 contains further technical details and remarks, including a discussion of the importance of the bounded shrinkage assumption and results for the posterior mean estimator.

1.8. *Software.* The methods described in [33] and in the current paper have been implemented as the EbayesThresh contributed package within the R statistical language [45]. The package and documentation can be installed from the CRAN archive accessible from http://www.R-project.org. Additional description and implementational details are available in [34]. For a MATLAB implementation, see [6].

**2. Mixture priors and details of calculations.** In this section we discuss general aspects of the priors used in our procedure, and then review some theory for the single sequence case. Throughout, we use $c$ to denote generic strictly positive constants, not necessarily the same at each use, even within a single equation. When there is no confusion about the value of the prior weight $w$, it may be suppressed in our notation. We write $\Phi$ for the standard normal cumulative, and set $\tilde{\Phi} = 1 - \Phi$. It is assumed throughout that the model and the observed data are renormalized so that the noise variance $\sigma_E^2 = 1$.

2.1. *Priors with heavy tails.* Particular heavy-tailed densities that we shall consider for the nonzero part of the prior distribution are the Laplace density with scale parameter $a > 0$,

$$\gamma_a(u) = \tfrac{1}{2}a\exp(-a|u|),$$

and the mixture density given by

(12) $\qquad \mu|\Theta = \vartheta \sim N(0, \vartheta^{-1} - 1) \qquad \text{with } \Theta \sim \text{Beta}(\tfrac{1}{2}, 1).$

More explicitly, the latter density for $\mu$ has

(13) $\qquad \gamma(u) = (2\pi)^{-1/2}\{1 - |u|\tilde{\Phi}(|u|)/\varphi(u)\}$

and has tails that decay as $u^{-2}$, the same weight as those of the Cauchy distribution. For this reason we refer to the density (13) as the *quasi-Cauchy* density.

We shall mostly consider functions $\gamma$ that satisfy the following conditions:

1. The function $\gamma$ is a symmetric unimodal density satisfying the condition

(14) $\qquad \sup_{u>0}\left|\frac{d}{du}\log\gamma(u)\right| < \infty.$

2. The quantity $u^2\gamma(u)$ is bounded over all $u$.
3. For some $\kappa \in [1, 2]$,

$$y^{1-\kappa}\gamma(y)^{-1}\int_y^\infty \gamma(u)\,du$$

is bounded above and below away from zero for sufficiently large $y$.

The first of these conditions implies that the tails of $\gamma$ are exponential or heavier, while the second rules out tail behavior heavier than Cauchy. The third condition is a mild regularity condition. The conditions are satisfied if $\gamma$ is the Laplace or quasi-Cauchy function, but not if $\gamma$ is a normal density.

For the normal, Laplace and quasi-Cauchy priors, the posterior distribution of $\mu$, given an observed $Z$, and the marginal distribution of $Z$ are tractable, so that the choice of $w$ by marginal maximum likelihood, and the estimation of $\mu$ by posterior mean or median, can be performed in practice, as outlined in the following paragraphs. We begin by setting out generic calculations for the relevant quantities, and then give specific details for particular priors.

2.2. *Generic calculations.*

*Posterior mean.* In general, the posterior probability $w_{\text{post}}(z) = P(\mu \neq 0 | Z = z)$ will satisfy

(15) $$w_{\text{post}}(z) = wg(z)/\{wg(z) + (1-w)\varphi(z)\}.$$

Define
$$f_1(\mu|Z=z) = f(\mu|Z=z, \mu \neq 0),$$
so that the posterior density
$$f_{\text{post}}(\mu|Z=z) = (1 - w_{\text{post}})\delta_0(\mu) + w_{\text{post}} f_1(\mu|z).$$
Let $\mu_1(z)$ be the mean of the density $f_1(\cdot|z)$. The posterior mean $\tilde{\mu}(z; w)$ is then equal to $w_{\text{post}}(z)\mu_1(z)$.

*Posterior median.* To find the posterior median $\hat{\mu}(z;w)$ of $\mu$, given $Z = z$, let
$$\tilde{F}_1(\mu|z) = \int_\mu^\infty f_1(u|z)\,du.$$
If $z > 0$, we can find $\hat{\mu}(z,w)$ from the properties

(16) $$\begin{aligned}\hat{\mu}(z;w) &= 0 &&\text{if } w_{\text{post}}(z)\tilde{F}_1(0|z) \leq \tfrac{1}{2},\\ \tilde{F}_1(\hat{\mu}(z;w)|z) &= \{2w_{\text{post}}(z)\}^{-1} &&\text{otherwise.}\end{aligned}$$

Note that if $w_{\text{post}}(z) \leq \frac{1}{2}$, then the median is necessarily zero, and it is unnecessary to evaluate $\tilde{F}_1(0|z)$. If $z < 0$, we use the antisymmetry property $\hat{\mu}(-z, w) = -\hat{\mu}(z, w)$.

*Marginal maximum likelihood weight.* The explicit expression for the function $g$ facilitates the computation of the maximum marginal likelihood weight in the single sequence case. Define the score function $S(w) = \ell'(w)$, and define

(17) $$\beta(z, w) = \frac{g(z) - \varphi(z)}{(1-w)\varphi(z) + wg(z)}.$$

Then $\beta(z, w)$ is a decreasing function of $w$ for each $z$, and

(18) $$S(w) = \sum_{i=1}^n \beta(Z_i, w).$$

Letting $w_n$ be the weight that satisfies $t(w_n) = \sqrt{2 \log n}$, the estimated weight $\hat{w}$ maximizes $\ell(w)$ over $w$ in the range $[w_n, 1]$. It follows that, if $S$ has a zero in this range, then $S(\hat{w}) = 0$. Furthermore, the smoothness and monotonicity of $S(w)$ make it possible to find $\hat{w}$ by a binary search, or an even faster algorithm. The restriction on the range of $\hat{w}$ implies that the threshold $t(\hat{w}) \leq \sqrt{2 \log n}$.

*Shrinkage rules.* The posterior median and mean are examples of estimation rules that yield an estimate of $\mu$, given $Z = z$. In general, a family of estimation rules $\eta(z,t)$, defined for all $z$ and for $t > 0$, will be called a *thresholding rule* if and only if, for all $t > 0$, $\eta(z,t)$ is an antisymmetric and increasing function of $z$ on $(-\infty, \infty)$ and $\eta(z,t) = 0$ if and only if $|z| \leq t$. It will have the *bounded shrinkage* property if and only if

$$z - (t + b_0) \leq \eta(z,t) \leq z \qquad \text{for all } z > t \tag{19}$$

for some constant $b_0$ independent of $t$.

An immediate consequence of (19) is that $|z - \eta(z,t)| \leq t + b_0$ for all $z$ and $t$. For any given weight $w$, the posterior median will be a thresholding rule, with a threshold we denote by $t(w)$, and will have the bounded shrinkage property under condition (14). More general thresholding rules may have advantages in some cases. For example, the hard thresholding rule, with a suitably estimated threshold, may have computational advantages and may preserve peak heights better, but we have not investigated this aspect in detail. Indeed, the choice of shrinkage rule and the choice of threshold are somewhat separate issues. The former is problem dependent and this paper is devoted to the latter.

The posterior mean is not a thresholding rule, but has sufficient properties in common with the posterior median to allow similar theoretical results to be obtained, but under restrictions on the risk functions considered.

2.3. *Calculations for specific priors.* The calculations set out above show that the key quantities are the marginal density $g$, the mean function $\mu_1(z)$ and the tail conditional probability function $\tilde{F}_1$. If $\gamma$ is the $N(0, \tau^2)$ density, then $g$ will be the $N(0, 1 + \tau^2)$ density, and $\mu_1(z) = \lambda x$, where $\lambda = \tau^2/(1 + \tau^2)$. The function $\tilde{F}_1(\mu|x)$ will be the upper tail probability of the $N(\lambda x, \lambda)$ density.

For the Laplace distribution prior, we have

$$g(z) = \tfrac{1}{2} a \exp(\tfrac{1}{2} a^2) \{ e^{-az} \Phi(z - a) + e^{az} \tilde{\Phi}(z + a) \}$$

and

$$f_1(\mu|z)$$
$$= \begin{cases} e^{az} \varphi(\mu - z - a)/\{e^{-az}\Phi(z - a) + e^{az}\tilde{\Phi}(z + a)\}, & \text{if } \mu \leq 0, \\ e^{-az} \varphi(\mu - z + a)/\{e^{-az}\Phi(z - a) + e^{az}\tilde{\Phi}(z + a)\}, & \text{if } \mu > 0, \end{cases} \tag{20}$$

which is a weighted sum of truncated normal distributions. Hence, it can be shown that, for $z > 0$,

$$\mu_1(z) = z - \frac{a\{e^{-az}\Phi(z - a) - e^{az}\tilde{\Phi}(z + a)\}}{e^{-az}\Phi(z - a) + e^{az}\tilde{\Phi}(z + a)}. \tag{21}$$

For $\mu \geq 0$, under the Laplace prior, we have

$$\tilde{F}_1(\mu|z) = \frac{e^{-az}\tilde{\Phi}(\mu - z + a)}{e^{-az}\Phi(z-a) + e^{az}\tilde{\Phi}(z+a)}.$$

For the quasi-Cauchy distribution, we have

$$g(z) = (2\pi)^{-1/2} z^{-2}(1 - e^{-z^2/2})$$

and

$$\mu_1(z) = z(1 - e^{-z^2/2})^{-1} - 2z^{-1}.$$

After some manipulation,

$$\tilde{F}_1(\mu|z) = (1 - e^{-z^2/2})^{-1}\{\tilde{\Phi}(\mu - z) - z\varphi(\mu - z) + (\mu z - 1)e^{\mu z - z^2/2}\tilde{\Phi}(\mu)\}.$$

For the Laplace prior, the equation $\tilde{F}_1(\hat{\mu}(z;w)|z) = \{2w_{\text{post}}\}^{-1}$ in (16) can be solved explicitly for $\hat{\mu}(z;w)$, making use of the function $\Phi^{-1}$. In the case of the quasi-Cauchy prior, the equation has to be solved numerically.

**3. Some simulation results.** A simulation study was carried out for the regression models that are by now standard in the consideration of wavelet methods and are given in [18]. Simulations from each of the four models were carried out, for each of two noise levels. For "high noise," the ratio of the standard deviation of the noise to the standard deviation of the signal values is $\frac{1}{3}$. In the "low noise" case the ratio is $\frac{1}{7}$. This complements the simulations for the single sequence case reported in [33]. The S-PLUS code used to carry out the simulations is available from the authors' web sites, enabling the reader both to verify the results and to conduct further experiments if desired.

3.1. *Results for the translation-invariant wavelet transform.* In Table 1 various wavelet methods, all making use of the translation-invariant wavelet transform, are compared. For each model and noise level, 100 replications were generated. In each replication, the function was simulated at 1024 equally spaced points $t_i$. The same normal noise variables were used for each of the models and noise levels. The error reported for each method considered is

$$\sigma_E^{-2} \sum_{i=1}^{1024} \{\hat{f}(t_i) - f(t_i)\}^2,$$

where $\sigma_E^2$ is the noise variance in each case, and this explains why the results for "low noise" are apparently inferior to those for "high noise." The default choices of wavelet, boundary corrections and so on, given in the S-PLUS Wavelets function `waveshrink`, were used. For each realization, the noise variance is estimated using the median absolute deviations of the wavelet coefficients at the highest level. The default choice of boundary treatment is to use periodic boundary conditions, and such boundary conditions have to be used for current implementations of the translation-invariant wavelet transform. Detailed consideration of the use of the idea of the translation-invariant transform, in combination with boundary correction, is an interesting idea for future research.

For the Laplace prior $\gamma$, with both $w$ and the scale parameter $a$ estimated level-by-level by marginal maximum likelihood from the data, estimates were constructed using both the posterior median and the posterior mean. For the quasi-Cauchy prior, estimates using the posterior median were calculated. The posterior median for the mixed Gaussian prior was also calculated; as for the Laplace prior, both $w$ and the scale parameter were estimated from the data.

Three other methods based on the translation-invariant wavelet transform were considered: SURE applied to 4 and 6 levels of the transform, universal soft thresholding applied to 6 levels of the transform, and the false discovery rate approach with various values of the parameter $q$. Whenever the false discovery approach is used in the wavelet context, the method is applied separately at each level, a method derived from [2]. The same parameter $q$

TABLE 1
*Average over* 100 *replications of summed squared errors over* 1024 *points for various models and methods. All the wavelet-based estimators use the translation-invariant wavelet transform. The standard error of each of the entries is at most* 2% *of the value reported*

| Method | High noise | | | | Low noise | | | |
|---|---|---|---|---|---|---|---|---|
| | bmp | blk | dop | hea | bmp | blk | dop | hea |
| Laplace (median) | 171 | 176 | 93 | 41 | 212 | 164 | 109 | 57 |
| Quasi-Cauchy (median) | 177 | 185 | 97 | 40 | 221 | 169 | 115 | 56 |
| Gaussian (median) | 223 | 178 | 108 | 42 | 296 | 247 | 150 | 65 |
| Laplace (mean) | 181 | 182 | 100 | 45 | 214 | 175 | 115 | 62 |
| SURE (4 levels) | 243 | 205 | 140 | 73 | 299 | 255 | 181 | 95 |
| SURE (6 levels) | 237 | 199 | 123 | 45 | 296 | 252 | 167 | 71 |
| Univ soft (6 levels) | 701 | 417 | 229 | 67 | 997 | 749 | 386 | 110 |
| FDR ($q = 0.01$) | 170 | 198 | 97 | 43 | 223 | 164 | 109 | 56 |
| FDR ($q = 0.05$) | 169 | 173 | 93 | 39 | 223 | 163 | 110 | 53 |
| FDR ($q = 0.1$) | 177 | 168 | 93 | 39 | 235 | 174 | 116 | 53 |
| FDR ($q = 0.4$) | 264 | 212 | 127 | 50 | 353 | 273 | 181 | 72 |
| Spline | 1294 | 433 | 265 | 51 | 6417 | 1826 | 905 | 117 |
| Tukey | 545 | 330 | 286 | 246 | 1892 | 655 | 425 | 257 |

is used at each level, but the resulting estimated threshold may, of course, vary.

Comparisons are also included with two standard nonwavelet methods: cubic smoothing splines using GCV (`smooth.spline` in S-PLUS) and Tukey's 4(3RSR)2H method, running medians with twicing, the default S-PLUS `smooth`.

The standard error of each of the entries in the table is at most 2% of the value reported, so the values are correct to about 2 significant figures. The two standard nonwavelet methods both perform badly. Not surprisingly, given that it is specifically designed for smooth functions, the smoothing spline method fails disastrously on discontinuous and spiky signals. Neither method is good at separating signal from noise in the low noise case. The Tukey method is, to some extent, competitive with universal thresholding for the more inhomogeneous signals, but cannot adapt to the smoother behavior of the HeaviSine signal.

As for the methods based on the wavelet transform, the performance of the posterior mean estimator with the Laplace prior is consistently slightly worse than that of the posterior median. The universal thresholding method does not compare well, and SURE also gives noticeably worse performance than the Laplace and quasi-Cauchy empirical Bayes methods. The FDR method is competitive, provided the parameter is chosen appropriately. For these signals and sample size, $q = 0.05$ and $0.1$ give good performance, but the performance is worse in some cases if $q = 0.01$ and considerably worse if $q = 0.4$. We shall see in subsequent examples that the choice of this parameter is crucial to the performance of the FDR method, and that, in other situations, the relative performance of the FDR method is, in any case, not quite as good.

Within the translation-invariant wavelet transform, the observed coefficients are not independent. Benjamini and Yekutieli [8] propose a modification to the FDR method to take account of dependence between observations, replacing $q$ by $q/\sum_{k=1}^{M} k^{-1}$, where $M$ is the number of parameters under consideration. In the translation-invariant wavelet transform, the number of coefficients at each level is equal to the number of original observations, 1024 in the simulation example considered, so the correction factor is $\sum_{1}^{1024} k^{-1} \approx 7.5$. Therefore, the results reported for $q = 0.05$ would correspond to $q = 0.05 \times 7.5 = 0.375$ within the Benjamini–Yekutieli procedure. Since we are choosing the $q$ parameter arbitrarily in any case, this recalibration of the $q$ parameter does not affect our general conclusions. However, it does mean that the precise numerical value of $q = 0.05$ in the translation-invariant case cannot necessarily be translated directly to the standard discrete wavelet transform.

The mixed Gaussian prior model does not fit the theoretical assumptions of this paper and it can be seen that its performance is not as good as the

Table 2

*Difference in summed square errors between the methods indicated and the "Laplace (median)" method, measured in terms of the standard error of the difference estimated from* 100 *replications*

| Method | High noise | | | | Low noise | | | |
|---|---|---|---|---|---|---|---|---|
| | bmp | blk | dop | hea | bmp | blk | dop | hea |
| Quasi-Cauchy (median) | 14 | 16 | 10 | −2.9 | 16 | 12 | 15 | −2.6 |
| Laplace (mean) | 15 | 6 | 9 | 9 | 2.7 | 16 | 11 | 9 |
| SURE (6 levels) | 49 | 19 | 26 | 6 | 45 | 60 | 46 | 16 |
| Univ soft (6 levels) | 101 | 124 | 81 | 35 | 100 | 97 | 77 | 74 |
| FDR ($q = 0.01$) | −1.6 | 17 | 5 | 5 | 13 | −1.0 | 0.8 | −0.9 |
| FDR ($q = 0.05$) | −3.0 | −2.4 | −1.1 | −5 | 12 | −1.7 | 2.7 | −9 |
| FDR ($q = 0.1$) | 6 | −6 | −0.3 | −6 | 15 | 8 | 10 | −9 |
| FDR ($q = 0.4$) | 27 | 12 | 15 | 7 | 36 | 30 | 27 | 7 |

heavy-tailed priors. It is clear that the tail requirements on $\gamma$ have some bearing on the performance of the empirical Bayes approach. More detailed investigation of this issue would be an interesting topic for further research.

Because the same noise values are used for each model, there is correlation between the various values in Table 1. Comparisons of methods with the Laplace (median) method on a paired-sample basis are given in Table 2. It can be seen that the empirical Bayes method with the Laplace prior using the posterior median decisively outperfoms the other methods, except for the HeaviSine function, where the quasi-Cauchy prior performs very slightly better, but there is little to choose between the Laplace and quasi-Cauchy priors. Of the four FDR methods, the inferior performance for $q = 0.01$ and 0.4 is significant. For $q = 0.05$ and 0.1, the results are more equivocal, but the cases for which the FDR method underperforms are the ones with the most significant difference. Some further comparisons between these FDR methods and the empirical Bayes methods will be made below.

3.2. *Results for the standard discrete wavelet transform.* In order to evaluate the advantage of the translation-invariant transform, the same simulated data were also smoothed using methods based on the standard transform. The results are shown in Table 3. Additional comparisons are included, with the two block thresholding methods considered by Cai and Silverman [10], and with the QL method of Efromovich [23]. The block thresholding methods choose thresholds by reference to information from neighboring coefficients within the transform. In the case of NeighCoeff, only the two neighboring coefficients are used when considering a particular coefficient, while, for NeighBlock the data are processed in blocks and information is drawn from neighboring blocks. At coarse scales the QL method uses a thresholding rule with threshold equal to the standard deviation of the coefficients, while at finer levels the coefficients are thresholded at a threshold that increases up to the universal threshold as the level increases, but at the same time the proportion of coefficients allowed to be nonzero is also controlled, more stringently the higher the level.

TABLE 3
*Average over* 100 *replications of summed squared errors over* 1024 *points for various models and methods. In each case a standard wavelet transform was used. The two nonwavelet methods are not included, because they give the same results as in Table* 1. *For comparison, the results for the Laplace prior using the translation-invariant transform are repeated from Table* 1, *in italics*

| | High noise | | | | Low noise | | | |
|---|---|---|---|---|---|---|---|---|
| **Method** | bmp | blk | dop | hea | bmp | blk | dop | hea |
| *Laplace (median)* | | | | | | | | |
| *translation-invariant* | *171* | *176* | *93* | *41* | *212* | *164* | *109* | *57* |
| Laplace (median) | 278 | 245 | 147 | 53 | 338 | 311 | 204 | 76 |
| Quasi-Cauchy (median) | 277 | 252 | 150 | 54 | 324 | 301 | 200 | 73 |
| Gaussian (median) | 328 | 252 | 158 | 56 | 400 | 361 | 241 | 87 |
| Laplace (mean) | 257 | 228 | 140 | 57 | 304 | 278 | 190 | 79 |
| NeighBlock | 462 | 406 | 148 | 67 | 436 | 485 | 207 | 125 |
| NeighCoeff | 324 | 320 | 145 | 60 | 316 | 345 | 207 | 91 |
| QL | 359 | 310 | 175 | 58 | 411 | 366 | 243 | 82 |
| SURE (4 levels) | 317 | 248 | 183 | 97 | 393 | 331 | 247 | 117 |
| SURE (6 levels) | 312 | 247 | 167 | 69 | 399 | 339 | 235 | 94 |
| Univ soft (6 levels) | 937 | 484 | 277 | 76 | 1444 | 931 | 534 | 121 |
| FDR ($q = 0.01$) | 331 | 307 | 169 | 60 | 387 | 382 | 231 | 83 |
| FDR ($q = 0.05$) | 299 | 278 | 163 | 57 | 347 | 334 | 216 | 78 |
| FDR ($q = 0.1$) | 301 | 271 | 162 | 60 | 356 | 330 | 221 | 81 |
| FDR ($q = 0.4$) | 395 | 333 | 221 | 97 | 477 | 420 | 310 | 130 |

Several interesting conclusions can be drawn from this table. In this case, the posterior mean generally yields superior estimates to the posterior median. The NeighCoeff method is the better of the two block thresholding methods, but generally underperforms the Laplace prior/posterior mean method. The QL method performs well for the HeaviSine signal, but for the others is not so competitive. In this context, the relative performance of the FDR method is not as good as previously, but the importance of choosing the parameter $q$ appropriately remains. In general, it is clear how important is the use of a translation-invariant transform. The empirical Bayes method with a Gaussian prior was also tried in this context, and the results were, again, somewhat inferior to those for the heavy-tailed priors.

We can use Tables 1 and 3 to give another measure of performance. Let $r_{jk}$ denote the value in cell $(j,k)$ of the table, the error measure of method $j$ applied in case $k$. Then define the overall performance of method $j$ by $R(j) = \min_k(\min_\ell r_{\ell k}/r_{jk})$. The ratio $\min_\ell r_{\ell k}/r_{jk}$ quantifies the relative performance of method $j$ on case $k$, by comparing it with the best method for that case. The *minimum efficiency score* $R(j)$ then gives the loss of efficiency of estimator $j$ on the most challenging case. For the translation-invariant transform, the Laplace (median) case has a minimum efficiency score of 93%, while the FDR method with $q = 0.05$ scores 95%. The quasi-Cauchy method scores 91% and the FDR with $q = 0.1$ scores 90%.

However, if we turn to the standard transform, the results are more decisive, with scores of about 90% for the empirical Bayes Laplace and quasi-Cauchy median methods, but only 82% for the FDR with $q = 0.05$ and 84% for FDR with $q = 0.1$. It should be noted that the scores of around 90% for the empirical Bayes methods are only because the empirical Bayes method that is very best varies slightly between cases considered. But to be specific, the Laplace (median) method consistently outperforms all the FDR methods.

**4. Comparisons on illustrative data sets.** In this section the simulations are complemented by the consideration of three illustrative examples drawn from practical applications. Taking account of both the simulations and the practical comparisons, the empirical Bayes method, using the Laplace prior and the posterior median estimate, is fully automatic and, on each of the simulation studies considered as a whole, and on the practical illustrations, performs either best or nearly as well as the best method in each setting. The FDR method with $q = 0.05$ is slightly superior on the first simulation study, but at the expense of more substantial underperformance otherwise, at least on the cases we have considered.

4.1. *Inductance plethysmography data.* Our first practical comparison uses the inductance plethysmography data described in [39]. The data were collected by the Department of Anaesthesia, Bristol University, in an investigation of the breathing of patients after general anaesthesia. For further details, and the data themselves, see the help page for the ipd data in the WaveThresh package [40].

Plots of the original data and the curve estimate obtained using the Laplace prior method are shown in Figure 1. The results for the Laplace and quasi-Cauchy priors are virtually identical, so only the Laplace results are reported in detail here. The aim of adaptive smoothing with data of this kind is to preserve features such as peak heights as far as possible, while eliminating spurious rapid variation elsewhere. Abramovich, Sapatinas and Silverman [4] found that their BayesThresh method performed better in this regard than various other wavelet methods, but that for best results a subjective adjustment of their parameter $\alpha$ from $\alpha = 0.5$ to $\alpha = 2$ gave preferable results. The MML approach gave virtually the same results whether the quasi-Cauchy or Laplace prior is used.

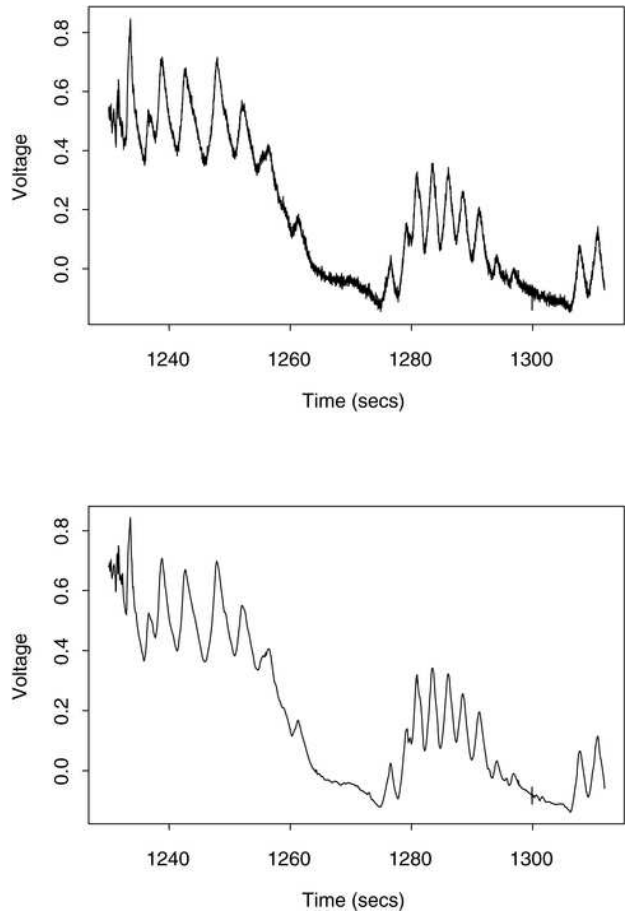

Fig. 1. Top panel: *the inductance plethysmography data.* Bottom panel: *the effect of smoothing the inductance plethysmography data with the Laplace prior method.*

The efficacy of the various methods in preserving peak heights is most simply judged by the maximum of the various estimates, the height of the first peak in the curve. The standard BayesThresh method ($\alpha = 0.5$) yields a maximum of 0.836, while subjectively adjusting to $\alpha = 2$ gives 0.845. The empirical Bayes method gives 0.842. Overall, the empirical Bayes method gives results much closer to the adjusted BayesThresh; the maximum distance from the empirical Bayes curves to the adjusted BayesThresh curve is about one-third that from the original BayesThresh estimate. The efficacy of the various methods in dealing with the rapid variation near time 1300 can be best quantified by the range of the estimated functions over a small interval near this point. The standard BayesThresh method has a "glitch" of range 0.08, while, for both the adjusted BayesThresh and the empirical Bayes method, the corresponding figure is under 0.06, a substantial if not dramatic improvement.

The FDR method with various parameters was also applied. As in the simulations, the FDR approach is applied separately to each level, with the same parameter $q$ at each level. For all the FDR $q$ parameters considered, the maximum of the estimated curve is between 0.842 and 0.843, but the range of the estimated curve near time 1300 is around 0.075. Thus, FDR competes well with empirical Bayes on preserving the peak height, but at the cost of inferior treatment of presumably spurious variation elsewhere.

Another comparison between the various methods can be made by considering the threshold that they use at various levels of the transform. The threshold is not a full description of the procedure, especially in the BayesThresh and Laplace prior cases where there are two parameters in the prior, but the threshold is a useful univariate summary of a method of processing wavelet coefficients. Figure 2 gives the comparison for the various methods applied to these data. It can be seen that, at the top four levels, the empirical Bayes methods track the adjusted BayesThresh method quite closely. The standard BayesThresh uses very high thresholds, which may be the reason why it smooths out the peak height somewhat. At the coarser levels, the empirical Bayes methods automatically adjust to much lower thresholds, reflecting a way in which the signal is less sparse at these levels, and thus allowing variation at these scales to go through quite closely to the way it is observed. None of the FDR parameter choices gives the degree of adaptivity of threshold to level shown by the empirical Bayes methods.

To conclude the comparison between BayesThresh and the empirical Bayes method, the subjectively adjusted BayesThresh method already yielded very good results for these data, but the basic message of this discussion is that the empirical Bayes method yields results virtually as good as those of the best BayesThresh method, but without any need for subjective tinkering with the parameters. In addition, the use of maximum likelihood to estimate the prior parameters is a less ad hoc approach than the fitting method used by the BayesThresh approach.

4.2. *Ion channel data.* A comparison between empirical Bayes and SURE is provided by considering a segment of the ion channel data discussed, for example, by Johnstone and Silverman [30]. Because these are constructed data, the "true" signal is known. See Figure 3. The thresholds chosen by SURE (dashed line) are reasonable at the coarse scales 6, 7 and 8, but are too small at the fine scales 9 to 11 where the signal is sparse, show some instability in the way they vary, and lead to insufficient noise removal in the reconstruction. By contrast, the empirical Bayes threshold choices increase monotonically with scale in a reasonable manner. In particular, the universal thresholds at levels 9 to 11 are found automatically. Two reconstructions using the same EB thresholds are shown in panel (b): one using the posterior median shrinkage rule, and the other using the hard thresholding rule. The hard threshold choice tracks the true signal better. The choice of threshold shrinkage rule is problem dependent, and beyond the scope of this paper. It is somewhat separate from the issue of setting threshold values.

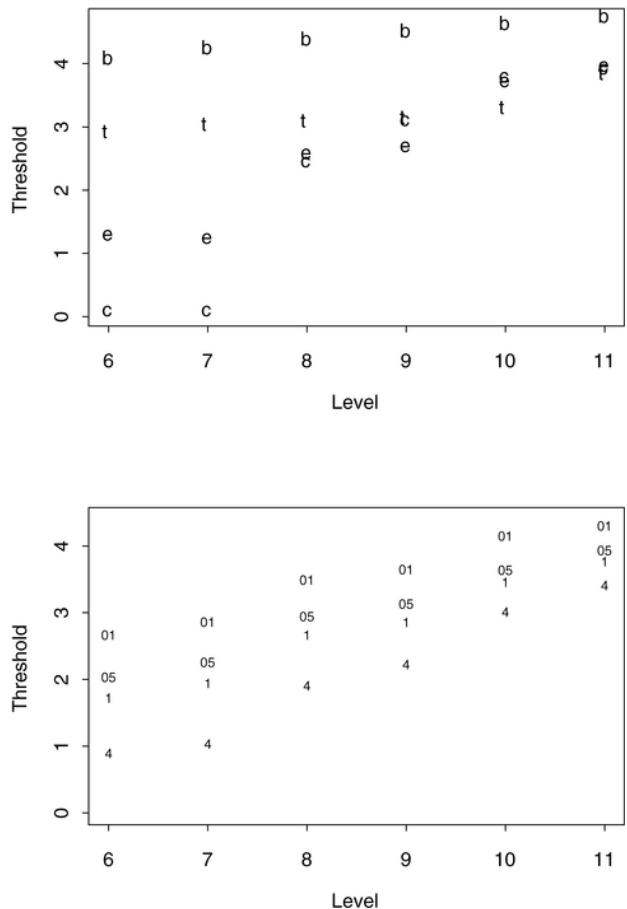

FIG. 2. *Thresholds chosen for the top six levels of the wavelet transform of the inductance plethysmography data by various methods.* Upper figure: e: *empirical Bayes, Laplace prior;* c: *empirical Bayes, quasi-Cauchy prior;* b: *BayesThresh;* t: *BayesThresh, subjectively tinkered, with* $\alpha = 2$. Lower figure: *False Discovery Rate method with parameters* $q = 0.01, 0.05, 0.1$ *and* $0.4$.

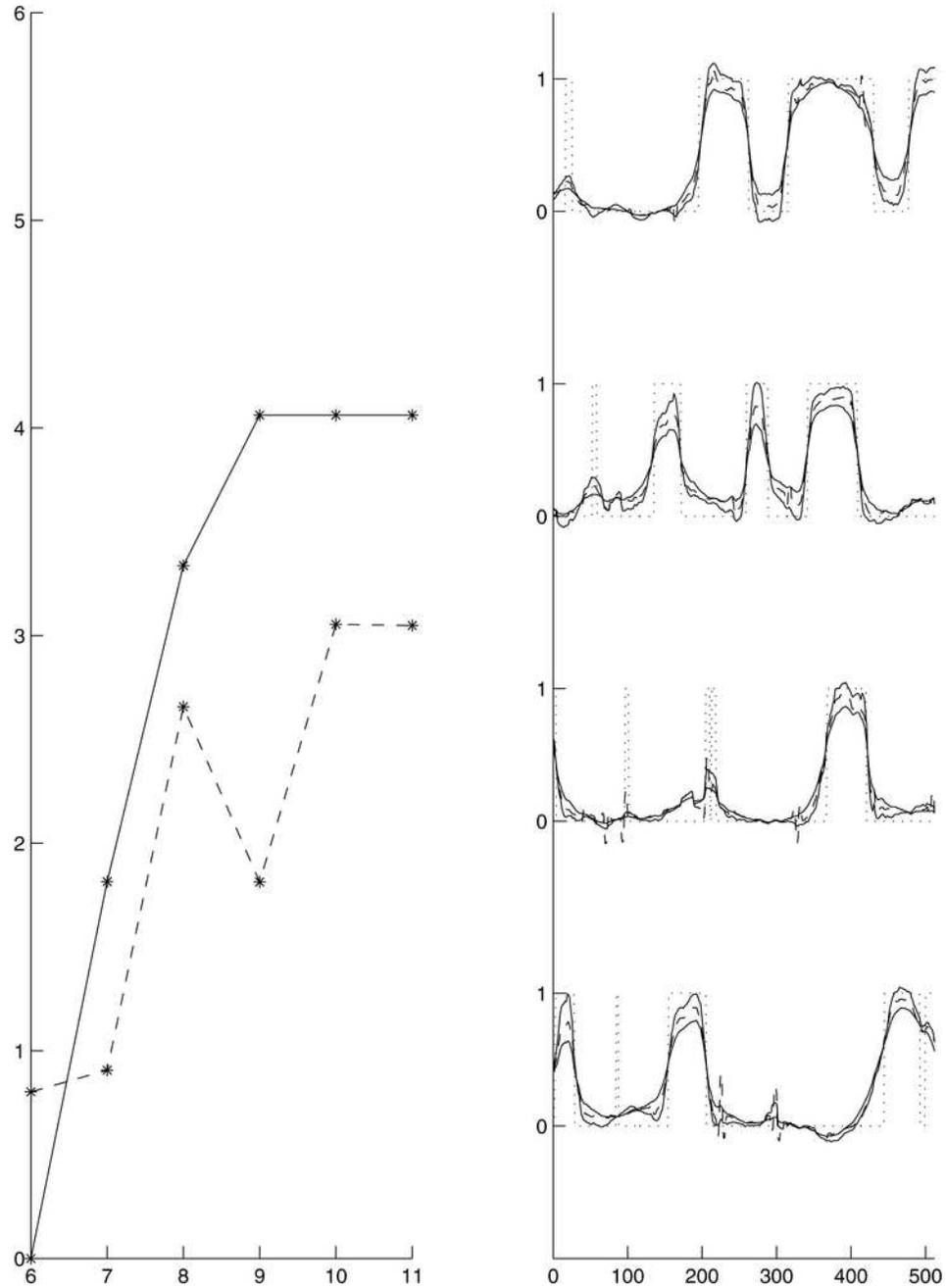

Fig. 3. Left panel: *Estimated threshold $\hat{t}$ plotted against level $j$; dashed line: SURE thresholds, solid line: EB thresholds.* Right panel: *Segment of the ion channel signal and three estimates. Both solid lines use EB-thresholds, but one uses a hard thresholding rule and tracks the true signal better, while the other uses posterior median shrinkage. The result of using SURE thresholds is plotted as the dashed line, and the dotted line gives the true signal.*

TABLE 4
*Percentage of errors in estimation of ion channel gating signal. The errors are the average over ten separate sequences of length* 4096 *drawn from the data provided by Eisenberg and Levis. The variances of the wavelet coefficients at each level were estimated separately*

| Decimated?           | N   | Y   |
|----------------------|-----|-----|
| Laplace (median)     | 2.4 | 3.0 |
| Quasi-Cauchy (median)| 2.7 | 3.5 |
| Laplace (mean)       | 2.3 | 2.6 |
| SURE (4 levels)      | 2.2 | 3.1 |
| SURE (6 levels)      | 2.3 | 3.2 |
| Univ soft (6 levels) | 6.0 | 7.5 |
| FDR ($q = 0.01$)     | 3.1 | 4.4 |
| FDR ($q = 0.05$)     | 2.8 | 3.9 |
| FDR ($q = 0.1$)      | 2.6 | 3.7 |
| FDR ($q = 0.4$)      | 2.3 | 3.6 |
| Spline               | 4.4 |     |
| Tukey                | 11  |     |
| AWS                  | 6.2 |     |
| Special              | 2.0 |     |

A systematic quantitative comparison is given in Table 4. For each method considered, ten sequences of length 4096 drawn from the original data were analyzed. The variances of the wavelet transform at the various levels were estimated by separate consideration, imitating the effect of using a sequence of observations with no signal to calibrate the method. For each method, the curve estimated by the smoothing method was then rounded off to the nearest of zero and one to give the final estimate. The figures given are the average percentage error over the ten sequences considered.

As an aside, we note that our theoretical results, of course, do not specifically include this zero–one loss of the estimate rounded to the nearer of zero or one. However, we do consider $L_q$ losses for $q$ near zero, which catch something of the flavor of discrete losses, in view of the fact that the limit as $q \to 0$ of the $q$th power of the $L_q$ norm is a zero–one loss.

Comparisons were made with the special-purpose method developed specifically for this problem by the originators of the data, and with standard smoothing methods, including the AWS method of Polzehl and Spokoiny [43]. The special-purpose method achieves an error rate of 2.0%; because of the specificity of this method, it is perhaps not surprising that it cannot be beaten by the more general-purpose methods we consider, but some of the translation-invariant wavelet methods come close. In this case the posterior mean slightly outperforms the posterior median, and other good methods are SURE and FDR with $q = 0.4$. If we use the parameter values $q = 0.05$ and 0.1 appropriate in our main simulation, then the results are inferior, underlining the need to tune the FDR parameter to the problem at hand.

4.3. *An image example.* Turning finally and briefly to images, Figure 4 shows the effect of applying empirical Bayes thresholds to a standard image with Gaussian noise added. The thresholds are estimated separately in each channel in each level. Nine realizations were generated, and the signal to noise ratio of the estimates (SNR $= 20\log_{10}(\|\hat{f} - f\|_2/\|f\|_2)$) calculated for both thresholding at $3\sigma_E$ and for the empirical Bayes thresholds. Smaller SNR corresponds to poorer estimation, though, of course, this quantitative measure does not necessarily correspond to visual perception of relative quality. The actual images shown correspond to the median of the nine examples, ordered by the increase in SNR between the $3\sigma_E$ threshold approach and the empirical Bayes approach.

For the example shown, the EB thresholds are displayed in the table below. They increase monotonically as the scale becomes finer and yield SNR $= 33.83$. They are somewhat smaller in the vertical channel, as the signal is stronger there in the peppers image. Fixing the threshold at $3\sigma_E$ in all channels leads to small noise artifacts at fine scales (SNR $= 33.74$), while fixing the threshold at $\sigma_E\sqrt{2\log n}$ (not shown) leads to a marked increase in squared error (i.e., reduced SNR).

| **Channel/Level** | 3 | 4 | 5 | 6 | 7 |
|---|---|---|---|---|---|
| Horizontal | 0 | 1.1 | 2.3 | 3.2 | 4.4 |
| Vertical | 0 | 0 | 2.0 | 3.0 | 4.4 |
| Diagonal | 0 | 1.7 | 2.7 | 4.1 | 4.4 |

**5. Theoretical results.** We now turn to the theoretical investigation of the proposed empirical Bayes method for curve estimation using wavelets. In doing so we distinguish between various different models for observed wavelet coefficients and for the theoretical coefficients of interest. Suppose throughout that level $J$ is such that the sum of the lengths of all the coefficient vectors below level $J$ is equal to $N$.

5.1. *Models for the observed data.* In the *white noise* model, it is assumed that we have independent observations $Y_{jk} \sim N(\theta_{jk}, N^{-1})$ of the wavelet coefficients $\theta_{jk}$ themselves. Because of the orthonormality properties of the wavelet decomposition, observations of this kind would be obtained by carrying out a wavelet decomposition of the function $f(t) + N^{-1/2}\,dW(t)$, where $dW(t)$ is a white noise process. In our main theory, we only use the $Y_{jk}$ at levels $j < J$, setting coefficients at higher levels to be zero.

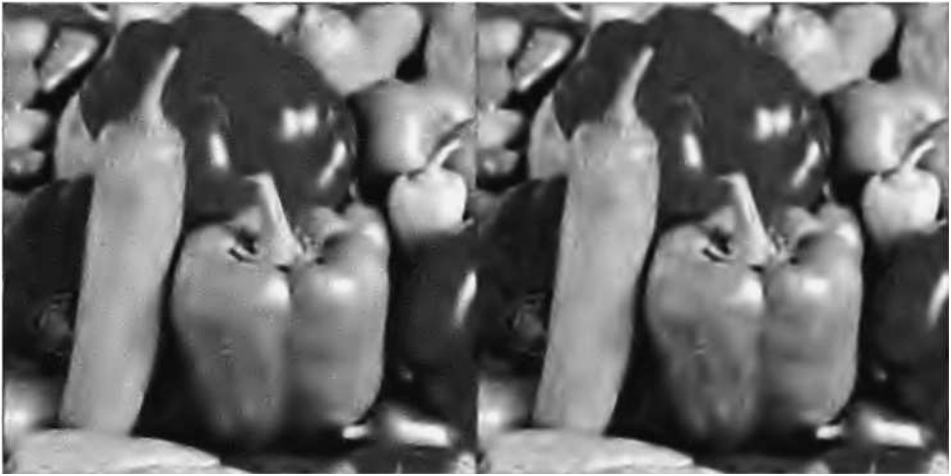

FIG. 4. *Translation invariant hard thresholding applied to a noisy version of the "peppers" image. For original image and noisy version see, for example, [36], Figure 10.6. Left panel: fixed threshold at $3\sigma_E$. Right panel: Level and channel dependent EB thresholds as shown in the table. The image obtained by fixed thresholds contains spurious high frequency effects that are largely obscured by the printing process. For a clearer comparison, the reader is recommended to view the images in the online version available from the authors' web sites.*

The other model of practical relevance is the *sampled data* model, where we assume that we have data $X_i = f(i/N) + \varepsilon_i$, where $\varepsilon_i$ are independent $N(0,1)$ random variables. Let $\tilde{\theta}$ be the discrete wavelet transform of the sequence $N^{-1/2} f(t_i)$, and $\tilde{Y}$ that of the sequence $N^{-1/2} X$, so that $\tilde{Y}_{jk} \sim N(\tilde{\theta}_{jk}, N^{-1})$. In much of the current statistical literature, the distinction between the white noise coefficients $Y_{jk}$ and the sampled-data coefficients $\tilde{Y}_{jk}$ is often glossed over, as is that between the function coefficients $\theta$ and the time-sampled coefficients $\tilde{\theta}$. The theoretical framework within which we work is, generally, to assume that

$$(22) \qquad \sum_k |\theta_{jk}|^p \leq C^p 2^{-apj} \qquad \text{for all } j,$$

corresponding to membership of the underlying function $f$ is a particular smoothness class. The first case we shall consider is where we observe $Y$ and estimate $\theta$.

The other cases all make use of the sampled-data coefficients $\tilde{Y}$. If we retain the constraint (22) on the underlying function, we can show that, provided the wavelet basis is chosen appropriately, the discretization involved in the sampled-data construction does not affect the order of magnitude of the accuracy of the estimates. This is the case whether we consider the estimates $\hat{\theta}(\tilde{Y})$ of the coefficients to be estimates of the wavelet coefficients $\theta$ of the function itself, or use the estimated coefficients to reconstruct an estimate of the sequence $f(i/N)$ via the discrete wavelet transform $\tilde{\theta}$. Unless we impose periodic boundary conditions, a key prerequisite for the consideration of the sampled data model is the development of appropriate boundary-corrected bases with corresponding preconditioning of the data near the boundaries, and we consider this aspect below.

A final model is the situation where it is the sequence of values $f(i/N)$ that is of primary interest, but we place the Besov array bounds on the discrete wavelet transform $\tilde{\theta}$ of this sequence rather than on the underlying function. We replace (22) by the constraint

$$\text{(23)} \qquad \sum_k |\tilde{\theta}_{jk}|^p \leq C^p 2^{-apj} \qquad \text{for all } j < J.$$

In this case we only require orthonormality of the discrete wavelet transform, but the condition (23) depends both on the function $f$ and on the particular $N$ under consideration. The asymptotic theorem should be thought of as a "triangular array" result, rather than a limiting result for a particular function $f$. The formalism of the proof is identical to the white noise case, except there is no need to consider terms for $j \geq J$ and this eliminates one of the error terms in the result.

5.2. *Array results under Besov body constraints.* Suppose that $\theta_{jk}$ is a coefficient array, defined for $j = L-1, L, L+1, \ldots$ and $0 \leq k < K_j$, for $K_j$ satisfying $2^{j-1} \leq K_j \leq 2^j$ for $j \geq L$ and $2^{L-1} \leq K_{L-1} \leq 2^L$. Let $N = \sum_{L-1 \leq j < J} K_j$ for integers $J$, and consider limits as $J \to \infty$. For given $J$, assume we have observations $Y_{jk} \sim N(\theta_{jk}, N^{-1}\sigma_E^2)$ for $j = L-1, L, \ldots, J-1$, $0 \leq k < K_j$. The variance $\sigma_E^2$ is assumed to be fixed and known, and without loss of generality we set $\sigma_E^2 = 1$.

Let $\theta_j$ denote the vector $(\theta_{jk} : 0 \leq k < K_j)$ and define $Y_j$ similarly. The vector $\theta_{L-1}$ is estimated by $Y_{L-1}$. For $L \leq j < J$, each vector $\theta_j$ is estimated separately by the empirical Bayes method described above. Set $\mu = N^{1/2}\theta_j$ and $Z = N^{1/2}Y_j$, and obtain an estimate of $\mu$ using a possibly modified threshold with parameter $A \geq 0$. If $A = 0$, then the threshold is not modified, while if $A > 0$, the threshold is as defined in (7). The threshold is that corresponding to the posterior median function, but provided this value of the threshold is used, the estimation can be carried out by any thresholding rule satisfying the bounded shrinkage property (19). We then set $\hat{\theta}_j = N^{-1/2}\hat{\mu}$. For $j \geq J$, finer scales than the observations assumed available, we set $\hat{\theta}_j = 0$.

The overall risk is defined to be

$$\text{(24)} \qquad R_{N,q,s}(\theta) = E\|\hat{\theta}_{L-1} - \theta_{L-1}\|_q^q + \sum_{j=L}^{\infty} 2^{sqj} E\|\hat{\theta}_j - \theta_j\|_q^q.$$

Under suitable conditions on the wavelet family, this norm dominates a $q$-norm on the $\sigma$th derivative of the original function if $s = \sigma + \frac{1}{2} - \frac{1}{q}$; see Section 5.6. The constant by which the contribution of the scaling coefficients is multiplied is somewhat arbitrary, and may be altered without affecting the overall method or results. We can now state the main result, which demonstrates that the empirical Bayes method attains the optimal rate of convergence of the mean $q$th-power error for all values of $q$ and $p$ down to 0.

The result also yields smoothness properties of the posterior estimate. It demonstrates, for values of $\sigma$ and $q$ satisfying the conditions of the theorem, that the coefficient array $\hat{\theta}$ has finite $b_{q,q}^{\sigma}$ norm and, hence, under suitable conditions on the wavelet has $\sigma$th derivative bounded in $q$-norm.

THEOREM 1. *Assume that $0 < p \leq \infty$ and $0 < q \leq 2$, and that $\alpha \geq \frac{1}{p}$. Suppose that the coefficient array $\theta$ falls in a sequence Besov ball $b_{p,\infty}^{\alpha}(C)$ so that*

(25) $$\|\theta_j\|_p \leq C 2^{-aj} \quad \text{for all } j,$$

*where $a = \alpha + \frac{1}{2} - \frac{1}{p} \geq \frac{1}{2}$. Let $s = \sigma + \frac{1}{2} - \frac{1}{q}$ and set*

$$r = (\alpha - \sigma)/(2\alpha + 1) \quad \text{and} \quad r' = (a - s)/2a.$$

*Assume that $\sigma \geq 0$ and that $\alpha - \sigma > \max(0, \frac{1}{p} - \frac{1}{q})$. Assume also that $sq \leq A$.*

*Then, for some quantity $c$ which does not depend on $C$ or $N$ (but may depend on $\alpha, p, \sigma, q$, as well as $\gamma$, $A$ and the wavelet family), the overall $q$-norm risk satisfies*

(26) $$R_{N,q,s}(\theta) \leq c\{\Lambda(C, N) + C^q N^{-r''q} + N^{-q/2} \log^\nu N\},$$

*where*

(27) $$\Lambda(C, N) = \begin{cases} C^{(1-2r)q} N^{-rq}, & \text{if } ap > sq, \\ C^{(1-2r')q} N^{-r'q} \log^{r'q+1} N, & \text{if } ap = sq, \\ C^{(1-2r')q} N^{-r'q} \log^{r'q} N, & \text{if } ap < sq, \end{cases}$$

$r'' = \alpha - \sigma - (\frac{1}{p} - \frac{1}{q})_+$ *and* $0 \leq \nu \leq 4$.

REMARKS. If $q \leq p$, then necessarily $ap > sq$ since $a > s$. However, if $q > p$, then the three cases in (27) correspond, respectively, to the "regular," "critical" and "logarithmic" zones described in [22].

Note first that, by elementary manipulations,

$$r' - r = \frac{ap - sq}{apq(2\alpha + 1)},$$

so the cases in (27) could equally be specified in terms of the relative values of $r'$ and $r$. Also, $a - s = \alpha - \sigma - \frac{1}{p} + \frac{1}{q} > 0$, so $r' > 0$ and $r'' = \min\{\alpha - \sigma, a - s\}$.

The condition $sq \leq A$ will be satisfied for all $q$ in $(0, 2]$ if $A \geq 2\sigma$. A particular situation in which this will hold is the "standard" case $A = 0$ and $\sigma = 0$.

The rates in (27) agree with the lower bounds to the minimax rates derived in Theorem 1 of [22], and so the first term of (26) is a constant multiple of the minimax dependence of the risk on the number of observations $N$ subject to the Besov body constraints. For fixed $C$ the other terms are of smaller order. The same rates arise in [17], which demonstrates that suitable estimators, dependent on $\alpha$, attain these rates for $q = 2$.

First consider the $N^{-r''q}$ term. Using the conditions $a \geq \frac{1}{2}$ and $\alpha > \sigma \geq 0$, we have $a - s = 2ar' \geq r'$ and $\alpha - \sigma = (2\alpha + 1)r > r$, so that $r'' \geq \min\{r, r'\}$. If $a > \frac{1}{2}$, then the inequality is strict and the $N^{-r''q}$ term will be of lower polynomial order than $\Lambda(C, N)$ in every case. If $a = \frac{1}{2}$ and $r' \leq r$, we will have $r'' = r'$, but, for fixed $C$, $\Lambda(C, N)$ will still dominate because of the logarithmic factor.

Since $r < \frac{1}{2}$, the $N^{-q/2} \log^\nu N$ term will always be of smaller order than $\Lambda(C, N)$. This term shows that, even if the Besov space constant $C$ is allowed to decrease as $N$ increases, or is zero, we have not shown that the risk can be reduced below a term of size $N^{-q/2}$, with an additional logarithmic term in certain cases. The exact definition of $\nu$ is

$$(28) \qquad \nu = \begin{cases} 0, & \text{if } sq < A, \\ (q+1)/2, & \text{if } sq = A > 0, \\ 3 + (q - p \wedge 2)/2, & \text{if } sq = A = 0. \end{cases}$$

*Truncating risk at fine scales.* Consider the estimation of $\tilde{\theta}$ from the transform $\tilde{Y}$, subject to the discretized constraints (23). In this case there is no need to consider levels $j \geq J$ in the risk, and the condition $a \geq \frac{1}{2}$, equivalent to $\alpha \geq 1/p$, can be relaxed to $a > 0$, equivalent to $\alpha > \frac{1}{p} - \frac{1}{2}$. Define

$$(29) \qquad \tilde{R}_{N,q,s}(f) = E\|\hat{\theta}_{L-1}(\tilde{Y}) - \tilde{\theta}_{L-1}\|_q^q + \sum_{j=L}^{J-1} 2^{sqj} E\|\hat{\theta}_j(\tilde{Y}) - \tilde{\theta}_j\|_q^q.$$

We then have the simpler result

$$(30) \qquad \tilde{R}_{N,q,s}(f) \leq c\{\Lambda(C, N) + N^{-q/2} \log^\nu N\}.$$

Define $\hat{f}(i/N)$ to be the sequence obtained by the inverse discrete wavelet transform applied to $N^{1/2}\hat{\theta}(\tilde{Y})$. In the "standard" case $\sigma = A = 0$ and $q = 2$, the orthogonality of the wavelet transform allows us to deduce from (30) that, subject to the constraint (23),

$$N^{-1} \sum_{i=1}^{N} E\{\hat{f}(i/N) - f(i/N)\}^2$$
$$\leq c\{C^{2/(2\alpha+1)} N^{-2\alpha/(2\alpha+1)} + N^{-1}(\log N)^{4-(1/2)(p \wedge 2)}\},$$

which implies (9).

*White noise model when fine scale observations are available.* If we assume that we have data $Y_{jk}$ for all levels, not just for $j < J$, then we can again relax the lower bound $a \geq \frac{1}{2}$ to $a > 0$. For definiteness, estimate $\hat{\theta}_j$ from the data for levels up to $j = \bar{J}^2$, and set the estimate to zero for higher levels. Then we will have the result

(31) $\quad R_{N,q,s}(\theta) \leq c\{\Lambda(C,N) + C^q \exp(-c' \log^2 N) + N^{-q/2} \log^{2\nu} N\}$

for a suitable $c' > 0$. The second term in (31) decays faster than polynomial rate in $N$ for any fixed $C$.

The proof of Theorem 1, together with the minor modifications required to prove (30) and (31), is given in Section 6.2 below.

5.3. *Wavelets whose scaling functions have vanishing moments.* Turn now to the issue of developing theory for the sampled data case subject to retaining the constraints on the function $f$ itself. Crucial to our theory are wavelets constructed from a scaling function $\phi$ with vanishing moments of order $1, 2, \ldots, R-1$, and $R$ continuous derivatives, for some integer $R$. The corresponding mother wavelet $\psi$ is orthogonal to all polynomials of degree $R-1$ or less, and both $\phi$ and $\psi$ are supported on the interval $[-S+1, S]$ for some integer $S > R$. Coiflets, as discussed in Chapter 8.2 of [16], are an example of wavelets constructed to satisfy these properties. The zero moments of the scaling function are used to control the discretization error involved when mapping observations to scaling function coefficients at a fine scale. Note that many standard wavelet families have scaling functions with nearly vanishing moments of orders 1 and 2; an issue for future investigation is the tradeoff in finite samples between relaxing the condition of exactly vanishing moments and using wavelets of narrower support than coiflets.

Unless we are happy to restrict attention to periodic boundary conditions, it is necessary to modify the wavelets and scaling functions near the boundary, and, hence, the filters used in the corresponding discrete wavelet transform. A construction following Section 5 of [14] can be used to perform this modification, while maintaining orthonormality of the basis functions. We review the application of the construction; for fuller details and properties see [32].

REMARKS. 1. If the restriction to (boundary modified) coiflets is needed for our theory, why is inferior behavior not observed for other Daubechies wavelet families in practice? In fact, it follows from [26] that if one recenters a Daubechies scaling function $\phi$ at its mean $\tau = \int x\phi(x)\, dx$, then the second moment necessarily vanishes. Thus, up to a horizontal shift $\tau$, one obtains two vanishing moments "for free."

2. The approach to sampled data taken by Donoho and Johnstone [20] works for a broad class of orthonormal scaling functions, by a less direct construction relating white noise and sampled data models through multiscale Deslauriers–Dubuc interpolation.

The construction is based on boundary scaling functions $\phi_k^B$ for $k = -R, -R+1, \ldots, R-2, R-1$, and boundary wavelets $\psi_k^B$ for $k = -S+1, -S+2, \ldots, S-1, S-2$. The support of these functions is contained in $[0, 2S-2]$ for $k \geq 0$ and in $[-(2S-2), 0]$ for $k < 0$. We fix a coarse resolution level $L$ such that $2S < 2^L$. At every level $j \geq L$, there are $2^j - 2(S-R-1)$ scaling functions, defined by

$$\phi_{jk}(x) = 2^{j/2}\phi_k^B(2^j x), \qquad k \in 0:(R-1),$$
$$\phi_{jk}(x) = 2^{j/2}\phi(2^j x - k), \qquad k \in (S-1):(2^j - S),$$
$$\phi_{jk}(x) = 2^{j/2}\phi_{k-2^j}^B(2^j(x-1)), \qquad k \in (2^j - R):(2^j - 1),$$

and $2^j$ wavelets

$$\psi_{jk}(x) = 2^{j/2}\psi_k^B(2^j x), \qquad k \in 0:(S-2),$$
$$\psi_{jk}(x) = 2^{j/2}\psi(2^j x - k), \qquad k \in (S-1):(2^j - S),$$
$$\psi_{jk}(x) = 2^{j/2}\psi_{k-2^j}^B(2^j(x-1)), \qquad k \in (2^j - S + 1):(2^j - 1).$$

All these functions are supported within $[0, 1]$. There are no scaling functions defined for $R \leq k < S - 1$ or for $2^j - S < k < 2^j - R$, but there are no such gaps in the definition of the wavelets. The $S - 1$ wavelets at each end are boundary wavelets, and have the same smoothness, on $[0, 1]$, and vanishing moments as the original wavelets. The $2^j - 2S$ interior wavelets are not affected by the boundary construction, and depend only on the $2^J - 2S$ interior scaling functions at the finest scale. At the coarsest level $L$, there will be $2^L - 2(S - R - 1)$ scaling coefficients; denote by $\mathcal{K}_{L-1}$ the set of indices for which the scaling functions $\phi_{Lk}$ are defined. At every level $j \geq L$, define $\mathcal{K}_j^B$ to be the set of $k$ for which $\psi_{jk}$ is a scaled version of a boundary wavelet, and $\mathcal{K}_j^I$ to be the set of $k$ for which $\psi_{jk}$ is a scaled version of $\psi$ itself.

Given a function $f$ on $[0, 1]$, we can now define the wavelet expansion of $f$ by

$$(32) \qquad f = \sum_{k \in \mathcal{K}_{L-1}} \theta_{L-1,k}\phi_{L,k} + \sum_{j=L}^{\infty} \sum_{k=0}^{2^j - 1} \theta_{jk}\psi_{jk},$$

where

$$\theta_{L-1,k} = \int_0^1 f(t)\phi_{Lk}(t)\,dt \qquad \text{for } k \text{ in } \mathcal{K}_{L-1}$$

and

$$\theta_{jk} = \int_0^1 f(t)\psi_{jk}(t)\,dt \qquad \text{for } j \geq L \text{ and } 0 \leq k < 2^j.$$

Where there is a need to distinguish between the boundary and interior wavelet coefficients, we write $\theta^I$ for the coefficients with $j \geq L$ and $k \in \mathcal{K}_j^I$ and $\theta^B$ for the boundary coefficients, those with $j \geq L$ and $k \in \mathcal{K}_j^B$.

5.4. *Constructing wavelet coefficients from discrete data.* Suppose now that we are given a vector of observations or of values of a function. In order to map these to scaling function coefficients at a suitable level, it is necessary to construct appropriate preconditioning matrices. In this section we define these matrices and set out certain of their properties. For more details, see [32].

On the left boundary define the $R \times R$ matrix $W$ and the $(S-1) \times R$ matrix $U$ by

$$W_{k\ell} = \int_0^\infty x^\ell \phi_k^B(x)\,dx, \qquad k=1,2,\ldots,R; \ell=0,1,\ldots,R-1,$$

$$U_{j\ell} = j^\ell, \qquad j=1,2,\ldots,S; \ell=0,1,\ldots,R-1.$$

Because $U$ is of full rank, we can define $A^L$ to be an $R \times (S-1)$ matrix such that $A^L U = W$. Similarly, the matrix $A^R$ is constructed to satisfy $A^R \bar{U} = \bar{W}$, where

$$\bar{W}_{k\ell} = \int_{-\infty}^0 x^\ell \phi_{-k}^B(x)\,dx, \qquad k=1,2,\ldots,R; \ell=0,1,\ldots,R-1,$$

$$\bar{U}_{j\ell} = (-1)^\ell j^\ell, \qquad j=1,2,\ldots,S; \ell=0,1,\ldots,R-1.$$

Given a sequence $X_0, X_1, \ldots, X_{N-1}$ with $N = 2^J$, define the preconditioned sequence $P_J X$ by

$$(P_J X)_k = \sum_{i=0}^{S-2} A_{ki}^L X_i, \qquad k \in 0:(R-1),$$

$$(P_J X)_k = X_k, \qquad k \in (S-1):(N-S),$$

$$(P_J X)_k = \sum_{i=1}^{S-1} A_{N-k,i}^R X_{N-i}, \qquad k \in (N-R):(N-1).$$

If the $X_i$ are uncorrelated with variance 1, then the variance matrix of the first part of $P_J X$ is $A^L (A^L)'$, while that of the last part, with indices taken in reverse order, is $A^R (A^R)'$.

There is some freedom in the choice of $A^L$ and $A^R$. For example, they can be defined such that not quite all the original sequence is needed to evaluate the preconditioned sequence. Specifically, to eliminate dependence on the first or last $S - R - 1$ values of theسequence, let $U_1$ be the square invertible matrix consisting of the last $R$ rows of $U$, and let $A^L = [0_{R \times (S-R-1)} : W U_1^{-1}]$, and $A^R$ correspondingly.

If, on the other hand, we have all the values in the sequence, then we have more freedom to choose $A^L$ and $A^R$. A natural choice is $A^L = WU^+$ and $A^R = \bar{W}\bar{U}^+$, where the superscript $^+$ denotes the Moore–Penrose generalized inverse. These choices will minimize the traces of the matrices $A^L(A^L)'$ and $A^R(A^R)'$ and, hence, the sum of the variances of $P_J X$, if we suppose that the $X_i$ are uncorrelated with unit variance.

In general, under the same assumption on $X$, let $c_A$ be the maximum of the eigenvalues of $A^L(A^L)'$ and $A^R(A^R)'$. Let $\tilde{Y}$ be the boundary-corrected discrete wavelet transform of the sequence $N^{-1/2} P_J X$. Then the eigenvalues of the variance matrix of $P_J X$ will be bounded by $c_A$. Because of the orthogonality of the boundary-corrected discrete wavelet transform, the variance of the elements of $\tilde{Y}$ is bounded by $c_A N^{-1}$.

The array $\tilde{Y}^I$ of interior coefficients $(\tilde{Y}_{jk}: L \le j < J, S-1 \le k \le 2^j - S)$ will only depend on $X_i$ for $S - 1 \le i \le N - S$, in other words, those $X_i$ left unchanged by the preconditioning. Therefore, $\tilde{Y}^I$ will be an uncorrelated array of variables with variance $N^{-1}$.

The preconditioning also makes it possible to get very good approximations to the wavelet coefficients of a smooth function $g$ from a sequence $g_J$ of discrete values $g_{Ji} = g(i 2^{-J})$. Define the vector $S_J g$ to be the preconditioned sequence $P_J g_J$. For any smooth function $g$, each $2^{-J/2} S_{Jk} g$ is a good approximation to the scaling coefficient $\int g \phi_{Jk}$ of $g$ at level $J$. To be precise, provided $g$ is $R$ times continuous differentiable on $[0,1]$, we have for each $k$,

$$(33) \quad \left| S_{Jk} g - 2^{J/2} \int_0^1 g(t) \phi_{Jk}(t)\, dt \right| \le c 2^{-JR} \sup\{|g^{(R)}(x)| : x \in \mathrm{supp}(\phi_{Jk})\}.$$

The result (33) depends on the vanishing moment properties of the scaling functions and on the construction of the preconditioning matrices. For full details, see Proposition 3 of [32].

5.5. *The boundary-corrected empirical Bayes estimator.* In this section we set out a detailed definition of a boundary corrected version of the empirical Bayes estimator, and prove that it has attractive theoretical properties. Assume throughout that a boundary corrected basis is in use.

Assume that for $N = 2^J$ we have sufficient observations

$$X_i = f(i/N) + \varepsilon_i, \qquad \varepsilon_i \text{ independent } N(0,1)$$

to evaluate the preconditioned sequence $P_J Y$. Let $\tilde{Y}$ denote the boundary corrected discrete wavelet transform of $N^{-1/2} P_J X$.

Define the estimated coefficient array $\hat{\theta}$ as follows:

- Estimate the coarse scaling coefficients by their observed values, so set

$$\hat{\theta}_{L-1} = \tilde{Y}_{L-1}.$$

- Estimate the interior coefficients $\hat{\theta}^I$ by applying the empirical Bayes method level-by-level to the observed array $\tilde{Y}^I$.
- Threshold the boundary coefficients separately. At level $j$, use a hard threshold of $\tau_A(j/N)^{1/2}$, where $\tau_A^2 \geq 2(1+A)c_A \log 2$, so that for each $k \in \mathcal{K}_j^B$

$$\hat{\theta}_{jk} = \tilde{Y}_{jk} I[|\tilde{Y}_{jk}| > \tau_A(j/N)^{1/2}|].$$

- For unobserved levels $j \geq L$, set $\hat{\theta}_{jk} = 0$.

In our main theoretical discussion, we measure the risk of this estimate as an estimate of the wavelet expansion of the function itself by

$$(34) \quad R_{N,q,s}^*(\theta) = E\|\hat{\theta}_{L-1} - \theta_{L-1}\|_q^q + \sum_{j=L}^{\infty} 2^{sqj} E\|\hat{\theta}_j - \theta_j\|_q^q.$$

If we use $\hat{\theta}$ as an estimate of the discrete wavelet transform array $\hat{\theta}$, then the natural measure of accuracy is the risk $\tilde{R}$ as defined in (29). However, it should be noted that because of the preconditioning, the array $\hat{\theta}$ only specifies uniquely the values of the sequence $f(i/N)$ away from the boundaries.

The main result of this section demonstrates that the estimate has optimal-rate risk behavior over $q$ and $p$ down to zero.

THEOREM 2. *Assume that the scaling function $\phi$ and the mother wavelet $\psi$ have $R$ continuous derivatives and support $[-S+1, S]$ for some integer $S$, and that $\int x^m \phi(x)\, dx = 0$ for $m = 1, 2, \ldots, R-1$. Assume that the wavelets and scaling functions are modified by the boundary construction described above, and that the thresholding is carried out by a modified threshold method with $A \geq 0$. Assume that the available data and the construction of the estimator are as set out above.*

*Assume that $\alpha > \sigma$, $a > s$, $\alpha < R$, and $sq \leq A$. Assume either that $\alpha > \frac{1}{p}$ or that $\alpha = p = 1$. Assume that $0 < p \leq \infty$ and $0 < q < 2$. Let $r = (\alpha - \sigma)/(2\alpha + 1)$ and $r' = (a-s)/(2a)$. Let $\mathcal{F}(C)$ be the set of functions $f$ whose wavelet coefficients fall in $b_{p,\infty}^\alpha(C)$. Then there is a constant $c$ independent of $C$ such that, for suitable $r'''$ and suitable $\lambda$ and $\lambda'$,*

$$(35) \quad \sup_{f \in \mathcal{F}(C)} R_{N,q,s}^*(f) \leq c\{\Lambda(C,N) + C^q N^{-r'''q} \log^{\lambda'} N + N^{-q/2} \log^\lambda N\},$$

*where $\Lambda(C,N)$ is as defined in (27) in Theorem 1, and, for all fixed $C$, $N^{-r'''q} \log^{\lambda'} N$ is of smaller order than $\Lambda(C,N)$.*

The general correspondence between Besov sequence and function norms is discussed in Section 5.6. The case $\alpha = p = 1$ is included because the space defined by membership of the Besov sequence space $b^1_{1,\infty}(C)$ is well known to include the space of functions with appropriately bounded total variation.

If our only concern is for the accuracy of estimation of the array $\tilde{\theta}$, then we have the bound

$$
(36) \qquad \sup_{f \in \mathcal{F}(C)} \tilde{R}_{N,q,s}(f) \leq c\{\Lambda(C,N) + N^{-q/2} \log^\lambda N\}.
$$

5.6. *Besov array norms and function norms.* Our theory gives minimax risk bounds over functions whose array of wavelet coefficients fall in $b^\alpha_{p,\infty}(C)$ Besov sequence balls. Under appropriate assumptions on the wavelet basis, these Besov sequence norms on the wavelet coefficients are equivalent to the corresponding Besov function norms on the functions themselves. Relevant results for the specific case of boundary-corrected wavelets on a bounded interval are considered in detail in Appendix D of [29]. The equivalence will certainly hold for the wavelets with bounded support and vanishing moments up to order $R-1$, provided $p \geq 1$, $R \geq 2$, $R > \alpha$ and $\alpha > \max(0, \frac{1}{p} - \frac{1}{2})$. See also [19, 20] and further literature referenced there.

The $B^\alpha_{p,\infty}$ Besov function norm is not very easy to grasp intuitively, and for integers $m$ it is helpful to compare the Sobolev space $W^m_p$, which has norm $(\int |f|^p + |f^{(m)}|^p)^{1/p}$. Our minimax results hold over balls in $B^m_{p,\infty}$ Besov spaces; for $p \geq 1$, the standard result that $W^m_p$ is embedded in the space $B^m_{p,\infty}$ shows that the results will hold for minimax rates over balls in the Sobolev space $W^m_p$.

Turn now to the error measure. Suppose that $f$ is a function with wavelet expansion given as in (32) above. Given any $\sigma \geq 0$ and $0 < q \leq 2$, define

$$
r_{q,\sigma}(f) = \|\theta_{L-1}\|^q_q + \sum_{j=L}^\infty 2^{sqj} \|\theta_j\|^q_q,
$$

where $\sigma = s - \frac{1}{2} + \frac{1}{q}$. This corresponds to the risk measure (24) relative to which our theoretical bounds are obtained. We state and prove a proposition showing that the error norm dominates various function seminorms. It follows from the proposition that, for $\sigma \geq 0$ and for $q$ in $(0,2]$, the bounds on estimation error proved for the Besov body error measure $r_{q,\sigma}(f)$ hold a fortiori for error measured by the integrated $q$th power of the derivatives up to order $\sigma$, provided that the wavelet satisfies appropriate regularity conditions. Note that the lower bounds on $\sigma$ and $q$ are zero, rather than the larger bounds on $\alpha$ and $p$ required for full Besov norm equivalence.

PROPOSITION 1. *Suppose $0 \leq \sigma < R$. For any integer $r$ such that $0 \leq r \leq \sigma$, and $q$ in $(0, 2]$,*

$$\int_0^1 |f^{(r)}(t)|^q \, dt \leq c r_{q,\sigma}(f).$$

REMARK. The one-sided bound leaves open the possibility that the Sobolev norm might be much smaller than the Besov $B_{q,q}^\sigma$ norm on the right-hand side. However, for $1 < q \leq 2$ there is a reverse inequality,

$$c_1 \|f\|_{B_{q,2}^\sigma} \leq \|f\|_{W_q^\sigma},$$

in terms of a Besov norm on the space $B_{q,2}^\sigma$ (defined, e.g., in [47]). Since the minimax rate of convergence is the same for the $B_{q,2}^\sigma$ and the $B_{q,q}^\sigma$ error norms on regular and logarithmic zones [22], it can be said that we are capturing the situation for the Sobolev norm without too much loss.

5.7. *Comparisons.* In this section we compare the theoretical results for empirical Bayes estimators established in Theorems 1 and 2 in this paper with those known for some other existing thresholds.

The universal threshold $\sigma_E \sqrt{2 \log N}$ of [21] leads to rates of convergence that are suboptimal by logarithmic terms in the regular case and some critical cases; see the detailed discussion in [22], Section 12.1. For example, using the notation of (26) and (27), the bound on the rate in the regular case $ap > sq$ would be $\Lambda(C, N) = C^{(1-2r)q}((\log N)/N)^{rq}$. The reason that the universal threshold is suboptimal in this way is essentially that, in dense cases, thresholds should be set at a bounded (and small) number of standard deviations, rather than being of order $\sqrt{2 \log N}$.

For the SURE threshold recalled in Section 1.6, Donoho and Johnstone [19] and Johnstone [27] establish asymptotic optimality results under the special conditions of squared error loss ($q = 2$) for estimating the function ($\sigma = 0$) over Besov bodies with $p \geq 1$. Since the SURE estimate chooses thresholds to optimize an unbiased estimate of mean squared error, it is possible with these restrictions to obtain not only optimal rates, but also to show that the limiting MSE is minimax optimal even at the level of constants among threshold estimators. By the same token, it is less clear that one could expect better than optimal rates for other loss functions (say $q < 2$)—and even the rate issue remains to be formally investigated.

A more serious restriction of SURE is reflected in the constraint $p \geq 1$. As is discussed in [3] and [33], and illustrated in Figure 5, the SURE criterion $\hat{U}(t)$ in (11) is far from smooth in $t$. The asymptotic oracle inequality for SURE in Theorem 4 of [19] contains an error term of crude order $n^{-1/2}$, and this term has prevented any optimality conclusions from being drawn when $p < 1$. The instability shown in Figure 5 seems to derive from the derivative discontinuity created by the threshold zone: similar plots were obtained when applying the SURE criterion to the posterior median rule to estimate thresholds.

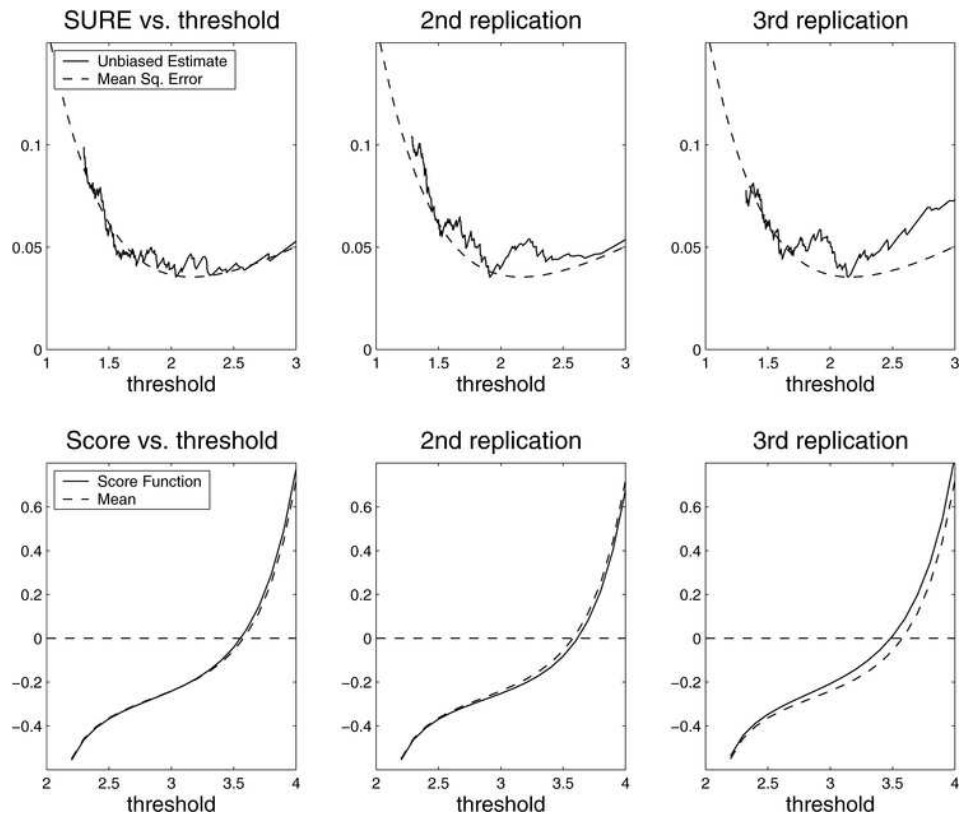

Fig. 5. *Instability of the SURE criterion compared to stability of the MML criterion. Three replications for $Z_k \stackrel{ind}{\sim} N(\mu_k, 1)$ with $\mu_k = 7$ for $k = 1:5$, zero for $k = 6:1000$. Top panels show the SURE criterion $\hat{U}(t)$ of (11) (solid) and its expectation (dashed). Bottom panels show the quasi-Cauchy score function (18) (solid) and its expectation (dashed) as a function of the (quasi-)threshold $t(w)$ defined by solving $\beta(t, 0) = 1/w$. (Johnstone and Silverman [33] has more on the quasi-threshold.)*

Related to this is another deficiency of the SURE threshold choice. While SURE adapts to squared error loss well on "dense" signals, the criterion does not reliably propose high thresholds for sparse signals. In order to obtain the theoretical results just cited, it was necessary to introduce a hybrid version of SURE containing a pretest for sparsity, which if detected, switched to a $\sqrt{2 \log n}$ threshold. Thus, the hybrid version creates a grossly discontinuous transition in thresholds which, while sufficient for the theoretical result, is unattractive in practice. Indeed, the simulations in [33] found the hybrid modification to be counterproductive in the examples considered.

Turn now to the levelwise FDR as applied in the wavelet context. While there has been extensive analysis of the exact adaptive minimax optimality of FDR in the sparse single sequence model [3] over $\ell_p$ balls and $\ell_q$ losses for $0 < p, q \leq 2$, there has been no published analysis of rates of convergence for a levelwise FDR estimate in the wavelet shrinkage setting. In unpublished work, IMJ combined the optimality properties of FDR for sparse signals with the advantages of SURE for dense signals using an improved pretest for sparsity. Adaptive optimality of rates of convergence was obtained under the conditions of Theorem 1 in the case $q = 2, \sigma = 0$. However, this work was abandoned in favor of the present empirical Bayes approach, due to the latter's smoother transition in threshold choice between sparse and dense regimes, reflected in better performance in actual examples.

Birgé and Massart [9] investigate a complexity-penalized model selection approach for Gaussian estimation and give a nonasymptotic risk bound for squared error loss (the case $q = 2, \sigma = 0$ here). When applied to our setting, their approach yields estimators that are minimax up to constants. The connection between the kind of "$2\log(n/k)$ per parameter" penalties used by Birgé and Massart [9] and Abramovich, Benjamini, Donoho and Johnstone [3] and FDR estimation is discussed further in the introduction to the latter paper. Finally, Paul [41] obtains optimal rates of convergence in certain inverse problems, again for $q = 2, \sigma = 0$, which in the direct estimation case would reduce to Theorem 1.

There has also been recent work on the optimality of Bayesian wavelet estimators based on mixture priors, though not from the adaptive estimation perspective: see, for example, [1, 42].

**6. Proofs of main results.** The proofs of the main theorems make use of results of Johnstone and Silverman [33] for the maximum marginal likelihood procedure in the single sequence case (3). In Section 6.1 we review these results, in a form recast to be useful for the multilevel problems raised in the current paper. In Section 6.2 we begin by giving an intuitive overview of the proof of our main result. This demonstrates the way that different kinds of error bounds are needed for different levels of the array; implicit in the proof is the way that the level-dependent empirical Bayes approach automatically adapts between these. After our intuitive discussion, the formal proof of Theorem 1 is given. The proof of Theorem 2 follows in Section 6.3; this makes use of approximation properties of the boundary-corrected bases and preconditioning operators defined in Sections 5.3 and 5.4 above.

6.1. *Results for the single sequence problem.* For a vector $\theta \in \mathbb{R}^n$, and $0 < q \leq 2$, suppose we have observations $Z \sim N_n(\theta, \varepsilon^2 I_n)$. Estimate $\theta$ by applying the marginal maximum likelihood approach to the data $\varepsilon^{-1}Z$ and then multiply the result by $\varepsilon$ to obtain an estimate of $\theta$. Assume that we are using a mixture prior (4) with $\gamma$ satisfying the assumptions set out in Section 2.1, and a family of thresholding rules with the bounded shrinkage property. We may use a modified thresholding method with $A > 0$, as defined in (7). By convention, let $A = 0$ denote the unmodified case. Define

$$(37) \qquad s_A^*(n) = \begin{cases} \log^{2+(q-p\wedge 2)/2} n, & \text{if } A = 0, \\ n^{-A}(\log n)^{(q-1)/2}, & \text{if } A > 0. \end{cases}$$

Then, by making appropriate substitutions into Theorem 2 of [33], we obtain the following result. The rates of convergence achieved by the leading terms in the various bounds are the minimax rates for the various parameter classes considered; see [33] for more details.

THEOREM 3. *Suppose that the above assumptions hold and that $0 < p \leq \infty$ and $0 < q \leq 2$. Then the estimate $\hat{\theta}$ satisfies the following risk bounds.*

(a) (Robustness.) *There exists a constant c such that*

$$(38) \qquad E\|\hat{\theta} - \theta\|_q^q \leq cn\varepsilon^q \qquad \text{for all } \theta.$$

(b) (Adaptivity for moderately sparse signals.) *There exist constants $c$ and $\eta_0$ such that, for sufficiently large $n$, provided $C_0^p < n\varepsilon^p \eta_0^p$, setting $\varepsilon_1 = \varepsilon\sqrt{\log(n\varepsilon^p C_0^{-p})}$,*

$$(39) \qquad \sup_{\|\theta\|_p \leq C_0} E\|\hat{\theta} - \theta\|_q^q \leq \begin{cases} c\{n^{1-(q/p)}C_0^q + \varepsilon^q s_A^*(n)\}, & \text{if } p \geq q, \\ c\{C_0^p \varepsilon_1^{q-p} + \varepsilon^q s_A^*(n)\}, & \text{if } p < q. \end{cases}$$

(c) (Adaptivity for very sparse signals.) *For $q > p > 0$ and $A > 0$, we also have, for sufficiently large $n$,*

$$(40) \qquad \sup_{\|\theta\|_p \leq C_0} E\|\hat{\theta} - \theta\|_q^q \leq c\{C_0^q + \varepsilon^q s_A^*(n)\} \qquad \text{for } C_0 < \varepsilon(\log n)^{1/2}.$$

*If $q \in (1,2]$, these results also hold if the estimation is carried out with the posterior mean function for the weight estimated by the marginal maximum likelihood procedure.*

In order to get an intuitive understanding of the way this result will be used in the multilevel setting, focus attention on the case $q > p$ and ignore the error term $\varepsilon^q s_A^*(n)$. The three results in Theorem 3 allow us to consider three zones of behavior of the underlying signal.

The first zone has large signal-to-noise ratio $C_0/\varepsilon > n^{1/p}\eta_0$. Here the best bound we have is a risk of order $n\varepsilon^q$, corresponding to the global risk of the maximum likelihood estimator $\hat{\theta}_{\text{MLE}}(Z) = Z$.

In the second zone, the signal-to-noise ratio is smaller and, since $q > p$, the risk can be substantially reduced by thresholding. A hard threshold rule with threshold $\varepsilon_1$ applied to the $Z$ will typically make an error in any individual coordinate of at most $\varepsilon_1$. A least favorable configuration satisfying the constraint $\|\theta\|_p \leq C_0$ would occur with $(C_0/\varepsilon_1)^p$ coordinates of size (a little less than) $\varepsilon_1$ each, and the rest being zero. This leads to a total error of order $(C_0/\varepsilon_1)^p \varepsilon_1^q$.

In the third zone, the region described by (40), the signal-to-noise ratio is so small that there is no benefit to attempting estimation at all, and $\hat{\theta}_0(Z) = 0$ is the natural estimator. This incurs risk $\|\theta\|_q^q \leq \|\theta\|_p^q \leq C_0^q$.

The discussion for $q \leq p$ is similar and simpler, involving only two zones, in the first of which the performance of the estimator is similar to that of $\hat{\theta}_{\mathrm{MLE}}$ and in the second to the zero estimator $\theta_0$. The impact of the theorem is that the empirical Bayes estimator adaptively achieves, roughly speaking, the best possible behavior whichever zone the signal actually falls in, without having to specify $p$ or $q$ or $C_0$ in advance.

We remark that, while Johnstone and Silverman [33] assumed $0 < p \leq 2$, we have subsequently checked that the results extend to $p \leq \infty$ and we use this broader range in this paper.

6.2. *Proof of Theorem* 1.

*Heuristic introduction.* Before the formal proof, we continue the intuitive discussion in order to give a heuristic explanation of where the rates of convergence and the phase change in the proof come from. In addition, we can gain an understanding of which kinds of estimators and bounds are needed (and, indeed, are imitated by our empirical Bayes method) at which levels of the transform. The discussion is inspired by the "modulus of continuity" point of view of Donoho, Johnstone, Kerkyacharian and Picard [22], but is adapted to the present setting. In the heuristic discussion, we ignore constants, error terms and so forth.

We apply the bounds of Theorem 3 level by level, with noise level $\varepsilon = N^{-1/2}$. The multiresolution index structure and Besov body constraints imply that, at level $j$, we have $n = K_j \asymp 2^j$ and $C_0 = 2^{-aj}C$. In the heuristic discussion approximate $\varepsilon_1$ by $(\log N/N)^{1/2}$ for simplicity, but in the actual proof this approximation is not used. With these substitutions, in the case $q > p$ the zones of Theorem 3 translate to

$$\sup_{\theta \in b_{p,\infty}^\alpha(C)} 2^{sqj} E\|\hat{\theta}_j - \theta_j\|_q^q$$

(41)

$$\approx \begin{cases} 2^{(sq+1)j} N^{-q/2}, & CN^{1/2} 2^{-(a+1/p)j} > 1, \\ 2^{-(ap-sq)j} C^p \left( \dfrac{N}{\log N} \right)^{-q/2}, & CN^{1/2} 2^{-(a+1/p)j} \leq 1 \\ & \text{and } C \left( \dfrac{N}{\log N} \right)^{1/2} 2^{-aj} > 1, \\ 2^{-(a-s)qj} C^q, & C \left( \dfrac{N}{\log N} \right)^{1/2} 2^{-aj} < 1. \end{cases}$$

The first zone corresponds to the coarsest scales; the transition to the middle zone occurs at scale $j_1$ defined by $2^{(\alpha+1/2)j_1} = CN^{1/2}$, and the third bound applies at scales above the finer index $j_2$ defined by $2^{aj_2} = C(N/\log N)^{1/2}$. The risk bounds increase geometrically as $j$ rises to $j_1$ and fall off geometrically as $j$ increases above $j_2$. The key to the behavior of the overall risk is $ap - sq$, because this determines the way the risk behaves in the zone between $j_1$ and $j_2$.

If $ap > sq$, then the least favorable index is $j_1$, and with geometric decay of risks away from this level, the rate is determined by $2^{(sq+1)j_1} N^{-q/2} = C^{1-2r} N^{-rq}$.

If $ap < sq$, the least favorable index is $j_2$ and the rate is given by $2^{(s-a)qj_2} C^q = C^{1-2r'} (\log N/N)^{r'q}$. Because of the extra logarithmic terms, this set of values of $(\alpha, \sigma, p, q)$ is referred to as the "logarithmic zone." Compare Figure 6.

If $ap = sq$, then each level between $j_1$ and $j_2$ contributes (in our heuristic approximation) an amount equal to the maximum in the case $ap < sq$. There are $O(\log N)$ such levels, leading to an extra $\log N$ factor.

We now give the formal proof, following this overall strategy.

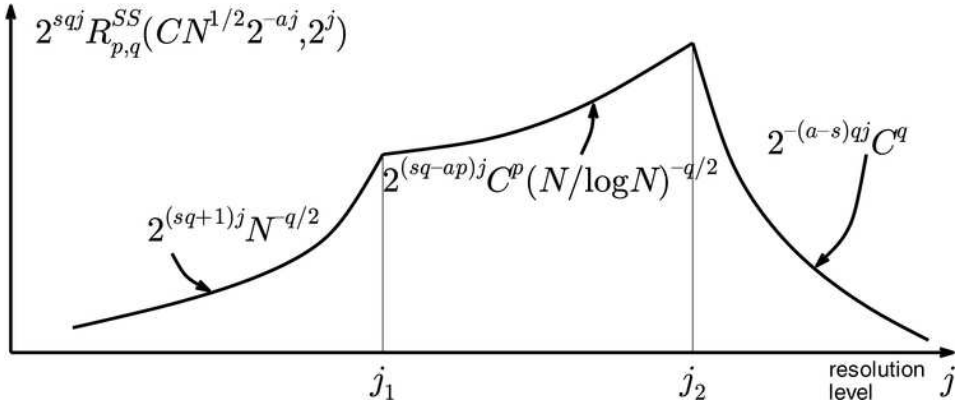

Fig. 6. *Schematic of risk contributions* (41) *by level $j$ in the "logarithmic" phase where $ap < sq$. In the setting $Z \sim N_n(\theta, \varepsilon^2 I_n)$ of Section 6.1, set $R_{p,q}^{SS}(C, n) = \inf_{\hat{\theta}} \sup_{\|\theta\|_p \leq C} E \|\hat{\theta} - \theta\|_q^q$.*

*Formal proof*: *division of the risk.* We split the sum in (24) into parts corresponding to the scaling coefficient, large-scale, fine-scale and very fine-scale parts of the risk. Define $n_0$ and $\eta_0$ so that (39) is satisfied. Let $j_1$ be the smallest integer for which $j \geq \max(L, \log_2 n_0)$ and $2^{1/p} \times 2^{-(a+1/p)(j+1)} C N^{1/2} < \eta_0$. We will then be able to apply the bound (39) at levels $j > j_1$. From the second property of $j_1$,

$$(42) \qquad 2^{-j_1} \leq c(C^2 N)^{-p/\{2(1+ap)\}} = c(C^2 N)^{-1/(2\alpha+1)}.$$

Also, $j_1$ must satisfy either $2^{j_1} \leq n_0 \vee 2^L$ or $2^{(\alpha+1/2)j_1} \leq 2^{1/p}\eta_0^{-1} C N^{1/2}$, so

$$(43) \qquad 2^{j_1} \leq c \max\{1, (C^2 N)^{1/(2\alpha+1)}\}.$$

Now write

$$(44) \qquad R_{N,q,s}(\theta) = E\|\hat{\theta}_{L-1} - \theta_{L-1}\|_q^q + R_{\mathrm{lo}} + R_{\mathrm{mid}} + R_{\mathrm{hi}},$$

where

$$R_{\mathrm{lo}} = \sum_{j=L}^{j_1 \wedge (J-1)} 2^{sqj} E\|\hat{\theta}_j - \theta_j\|^q,$$

$$R_{\mathrm{mid}} = \sum_{j_1 < j < J} 2^{sqj} E\|\hat{\theta}_j - \theta_j\|^q,$$

$$R_{\mathrm{hi}} = \sum_{j=J}^{\infty} 2^{sqj} \|\theta_j\|_q^q.$$

If $C^2 N \leq 1$, then in each case the last term in (26) will dominate the first. Therefore, we can assume throughout, without loss of generality, that $C > N^{-1/2}$, since there is no point in proving the result for any smaller values of $C$. It then follows from (42) and (43) that $2^{j_1} \asymp (C^2 N)^{1/(2\alpha+1)}$.

*Scaling function risk.* Let $b_q$ be the $q$th absolute moment of a standard normal random variable. Then

$$(45) \qquad E\|\hat{\theta}_{L-1} - \theta_{L-1}\|_q^q \leq K_{L-1} b_q N^{-q/2} = cN^{-q/2}.$$

Note that this holds regardless of the size of the scaling coefficients, and that the proof of Theorem 1 only uses the bound (25) for $j \geq L$.

*Risk at coarse scales.* To bound $R_{\mathrm{lo}}$, use the result (38) and the property $sq + 1 = q(\sigma + \frac{1}{2}) \geq \frac{1}{2}q > 0$ to give

$$(46) \qquad R_{\mathrm{lo}} \leq cN^{-q/2} \sum_{L \leq j \leq j_1} 2^{(sq+1)j} \leq cN^{-q/2} 2^{(1+sq)j_1}.$$

Therefore, since $(1+sq)/(2\alpha+1) = (\frac{1}{2}-r)q$,

(47) $$R_{\text{lo}} \leq cN^{-q/2}(C^2 N)^{(1/2-r)q} \leq cC^{(1-2r)q}N^{-rq}.$$

In the case $ap > sq$, this is exactly of the magnitude required in (26). If $ap \leq sq$, because $r \geq r'$ and $C^2 N > 1$, the bound (47) will be smaller than the first term in (26).

*Risk at fine scales.* As noted previously, by the definition of $j_1$, the bound (39) can be applied for the terms in the sum for $R_{\text{mid}}$. If $j_1 \geq J - 1$ there are no terms in the sum and $R_{\text{mid}} = 0$. Otherwise, set $\delta = (q-1)/2$ if $A > 0$ and $2 + (q - p \wedge 2)/2$ if $A = 0$ and define

(48) $$S_2 = N^{-q/2} \sum_{j_1 < j < J} 2^{sqj} s_A^*(K_j)$$
$$\leq N^{-q/2} \sum_{j_1 < j < J} 2^{-(A-sq)j} j^\delta \leq cN^{-q/2}(\log N)^\nu,$$

in every case, using the definition (28) of $\nu$.

Since at level $j$ we have $C_0 = C2^{-aj}$, considering the two terms in (39) now yields

(49) $$R_{\text{mid}} \leq c(S_1 + S_2),$$

where

$$S_1 = \begin{cases} C^q \sum_{j_1 < j < J} 2^{sqj} 2^{(1-q/p)j} 2^{-aqj}, & \text{if } p \geq q, \\ C^p N^{-(q-p)/2} \sum_{j_1 < j < J} 2^{-(ap-sq)j} \{\log(2^{(2\alpha+1)j} N^{-1} C^{-2})\}^{(q-p)/2}, & \\ & \text{if } q > p. \end{cases}$$

We consider four cases for $S_1$ and show that in every case $S_1 \leq c\Lambda(C, N)$. Combining with the bound (48) for $S_2$ allows us to conclude, in all cases, that

(50) $$R_{\text{mid}} \leq c\{\Lambda(C, N) + N^{-q/2}(\log N)^\nu\}.$$

*Case* 1a, $q \leq p$. In this case we necessarily have $ap > sq$ since $a > s$. If $p \geq q$, the exponent in the defining sum for $S_1$ is $(sq + 1 - q/p - aq)j = -(\alpha - \sigma)qj$, and, hence, is geometrically decreasing. Therefore,

(51) $$S_1 \leq cC^q 2^{-(\alpha-\sigma)qj_1} \leq cC^q(C^2 N)^{-rq} = cC^{(1-2r)q}N^{-rq}.$$

*Case* 1b, $q > p$ and $ap > sq$. The sum in $S_1$ is now geometrically decreasing apart from a log term, and so, using the property that $2^{(2\alpha+1)j_1} \asymp C^2 N$,

(52) $$S_1 \leq cC^p N^{-(q-p)/2} 2^{-(ap-sq)j_1} \{\log\{2^{(2\alpha+1)j_1} N^{-1} C^{-2}\}\}^{(q-p)/2}$$
$$\leq cC^p N^{-(q-p)/2} 2^{-(ap-sq)j_1} \leq cC^{(1-2r)q}N^{-rq}.$$

*Case* 2a, $q > p$ *and* $sq > ap$. In this case, necessarily $s > 0$ and, therefore, $A > 0$ since we require $sq \leq A$. Define $j_2$ by

$$2^{aj_2} = C(N/\log N)^{1/2}.$$

We now split the bounding sum for $S_1$ into the two zones $(j_1, j_2)$ and $[j_2, J)$. In the lower zone, since $C^2 N > 1$, we can bound

(53) $\quad \log(2^{(2\alpha+1)j} N^{-1} C^{-2}) \leq (2\alpha+1)j_2 \log 2 \leq c \log N \qquad \text{for } j \leq j_2.$

In the upper zone we use the bound (40).

The two sums obtained are set out in the following display. Since the terms in the first sum are geometrically increasing, and in the second geometrically decreasing, both sums are dominated by a multiple of their value at $j = j_2$:

(54) $\quad S_1 \leq cC^p (N^{-1} \log N)^{(q-p)/2} \sum_{j_1 < j < \min(j_2, J)} 2^{(sq-ap)j} + cC^q \sum_{j \geq j_2} 2^{sqj} 2^{-aqj}$

$\qquad \leq cC^p N^{-(q-p)/2} (\log N)^{(q-p)/2} 2^{(sq-ap)j_2} + cC^q 2^{-(a-s)j_2}$

(55) $\qquad \leq cC^{(1-2r')q} N^{-r'q} (\log N)^{r'q},$

after some algebra, substituting the definition of $j_2$.

*Case* 2b, $q > p$ *and* $ap = sq$. Argue as in Case 2a, but now the first sum in (54) is no longer geometric but is bounded by $j_2 \leq c \log N$. Carrying the extra $\log N$ factor through the argument yields

(56) $\qquad S_1 \leq cC^{(1-2r')q} (N/\log N)^{-r'q} \log N = c\Lambda(C, N).$

*Risk at very fine scales.* Define $\Delta = (\frac{1}{q} - \frac{1}{p})_+$, so that $r'' = a - s - \Delta$ and, whatever the relative values of $p$ and $q$, $\|\theta_j\|_q \leq 2^{\Delta j} \|\theta_j\|_p$ for each $j$. Then

(57)
$$R_{\mathrm{hi}} = \sum_{j=J}^{\infty} 2^{sqj} \|\theta_j\|_q^q \leq \sum_{j=J}^{\infty} 2^{sqj} 2^{\Delta qj} \|\theta_j\|_p^q$$

$$\leq C^q \sum_{j=J}^{\infty} 2^{-(a-s-\Delta)qj} = cC^q N^{-r''q}.$$

*Conclusion of proof and consideration of related results.* Combining the bounds (45), (47), (50) and (57) now completes the proof of Theorem 1.

To obtain (30) subject to the constraints (23), we follow exactly the same argument, noting that there are no coefficients to estimate for $j \geq J$, and so there is no term corresponding to (57).

To prove the result (31), where observations at all scales are available, modify the limits of summation where necessary. In the calculations for $R_{\mathrm{mid}}$, the sums are extended to $j = J^2$ where appropriate. None of the bounds for $S_1$ is affected, but for $S_2$ the calculation in (48) becomes

$$S_2 \leq N^{-q/2} \sum_{j_1 < j < J^2} 2^{-(A-sq)j} j^\delta \leq cN^{-q/2}(J^2)^\nu = cN^{-q/2} \log^{2\nu} N.$$

On the other hand, the sum in (57) now starts at $j = J^2$, leading to a bound

$$R_{\text{hi}} \leq cC^q 2^{-r''qJ^2} = cC^q \exp(-c \log^2 N).$$

Incorporating these two changes into the main argument leads to the result (31).

### 6.3. Proof of Theorem 2.

*Remarks and preliminaries.* In the estimation problem considered, the sequence $S_J f$ is the vector of expected values of the preconditioned data $P_J X$. Define $\tilde{\theta}$ to be the boundary-corrected discrete wavelet transform of $N^{-1/2} S_J f$. Our procedure uses what is essentially an estimate of $\tilde{\theta}$ to estimate the true coefficients $\theta$. The conditions of Theorem 2 allow the difference between these two arrays to be bounded; by Proposition 4 of [32] we have

(58) $\quad 2^{aj} \|\theta_j - \tilde{\theta}_j\|_p \leq cC 2^{-\bar{\alpha}(J-j)} \quad$ for all $j$ with $L - 1 \leq j < J$,

where $\bar{\alpha} = \alpha - (\frac{1}{p} - 1)_+ > \frac{1}{2}$. An immediate corollary is that, for some fixed constant $c$,

(59) $\quad \|\tilde{\theta}_j\|_p \leq cC 2^{-aj} \quad$ for $L - 1 \leq j < J$.

Therefore, the "discretized" coefficient array $\tilde{\theta}$ obeys (up to a constant) the same Besov sequence bounds as the "true" coefficient array $\theta$.

The precise value of the constant $r'''$ in the theorem is $\min\{a - s, \alpha - \sigma, \bar{\alpha}\} = \min\{r'', \bar{\alpha}\}$, with $\lambda' = 1$ if $\bar{\alpha} = \min(a - s, \alpha - \sigma)$ and 0 otherwise. We have already noted in the remarks following Theorem 1 that $N^{-r''q}$ is always of lower order than $\Lambda(C, N)$ for fixed $C$. The same is true of $N^{-r'''q} \log^{\lambda'} N$ since $\bar{\alpha} > \frac{1}{2}$ and the log term can only be present if $r''' > \frac{1}{2}$.

*Main component of error.* Use the convention that $I$ refers to the interior coefficients and $B$ to the boundary coefficients. The $\tilde{Y}_{jk}$ each have expected value $\tilde{\theta}_{jk}$, and for the interior coefficients are independent normals with variance 1. Because of the bound (59), we can argue exactly as in Theorem 1 to obtain

(60)
$$\sum_{j=L}^{J-1} 2^{sqj} E\|\hat{\theta}_j^I - \tilde{\theta}_j^I\|_q^q + \sum_{j>J} 2^{sqj} \|\theta_j\|_q^q$$
$$\leq c\{\Lambda(C, N) + C^q N^{-r''q} + N^{-q/2} \log^\nu N\},$$

where $\nu$ is as defined in (28).

Equation (60) gives the main part of the risk bound in Theorem 2, and the remainder of the proof consists in controlling all the other contributing errors.

*Coarse scale error.* Consider first the coarse level scaling coefficients $\theta_{L-1}$. Since the variance of each element of $\hat{\theta}_{L-1} = \tilde{Y}_{L-1}$ is bounded by $N^{-1} c_A$ we have

$$(61) \qquad E\|\hat{\theta}_{L-1} - \tilde{\theta}_{L-1}\|_q^q \leq cN^{-q/2} c_A^{q/2} 2^L \leq cN^{-q/2}.$$

*Boundary coefficients.* The contribution of the estimates of the boundary coefficients is considered in the following proposition.

PROPOSITION 2. *Under the assumptions of Theorem 2, uniformly over $\mathcal{F}(C)$ as $J \to \infty$,*

$$(62) \qquad \sum_{j=L}^{J-1} 2^{sqj} E\|\hat{\theta}_j^B - \tilde{\theta}_j^B\|_q^q \leq c\{C^{(1-2r')q}(N/\log N)^{-r'q} + N^{-q/2}\}.$$

PROOF. Define $R_{\text{boundary}}$ to be the sum on the left-hand side of (62). The array $\theta^B$ has the same number of coefficients at every level $j \geq L$, and the elements of the array $\tilde{Y}^B$ are normally distributed with expected values $\tilde{\theta}^B$, and variances bounded by $c_A N^{-1}$. We obtain the $\hat{\theta}_{jk}^B$ by individually hard thresholding the $\tilde{Y}_{jk}^B$ with threshold $\tau_A (j/N)^{1/2}$, so by standard properties of $q$-norms and thresholding,

$$(63) \qquad \begin{aligned} E|\hat{\theta}_{jk} - \tilde{\theta}_{jk}|^q &\leq c(E|\hat{\theta}_{jk} - \tilde{Y}_{jk}|^q + E|\tilde{Y}_{jk} - \tilde{\theta}_{jk}|^q) \\ &\leq c(j^{q/2} 2^{-qJ/2} + 2^{-qJ/2}) \\ &\leq cj^{q/2} 2^{-qJ/2}. \end{aligned}$$

If $s < 0$, use the bound (63) to give

$$(64) \qquad R_{\text{boundary}} \leq c 2^{-Jq/2} \sum_{j=L}^{J-1} j^{q/2} 2^{sqj} \leq cN^{-q/2}.$$

For $s \geq 0$, define $j_2$ by $2^{aj_2} = C(N/\log N)^{1/2}$ and split the risk into two parts. Letting $R_{\text{lo}}$ be the risk for boundary coefficients at levels below $\min(j_2, J)$, arguing as in (64),

$$(65) \qquad \begin{aligned} R_{\text{lo}} &\leq cN^{-q/2} \sum_{j=L}^{j_2 \wedge (J-1)} j^{q/2} 2^{sqj} \\ &\leq c(N/\log N)^{-q/2} 2^{sqj_2} \\ &\leq cC^{(1-2r')q}(N/\log N)^{-r'q}. \end{aligned}$$

Let $R_{\text{hi}}$ be the contribution to $R_{\text{boundary}}$ from levels $j \geq j_2$. For $j \geq j_2$, by Proposition 1 of [33], taking account of the bound $c_A$ on the variance of the $\tilde{Y}_{jk}$,

$$
\begin{aligned}
E|\hat{\theta}_{jk} - \tilde{\theta}_{jk}|^q &\leq c\{|\tilde{\theta}_{jk}|^q + c_A^{-1/2} j^{(q-1)/2} N^{-q/2} \varphi(c_A^{-1/2} \tau_A j^{1/2})\} \\
&\leq c\{|\tilde{\theta}_{jk}|^q + j^{(q-1)/2} N^{-q/2} 2^{-(1+A)j}\}.
\end{aligned}
\tag{66}
$$

Substituting the bound (66) and using the property $sq \leq A$ gives

$$
R_{\text{hi}} \leq c \sum_{j_2 \leq j < J} 2^{jsq} \|\tilde{\theta}_j^B\|_q^q + cN^{-q/2} \sum_{j \geq j_2} j^{(q-1)/2} 2^{-j}.
\tag{67}
$$

Since the vector $\tilde{\theta}_j^B$ is the same length for all $j$, for some constant $c$ independent of $j$ we have $\|\tilde{\theta}_j^B\|_q \leq c\|\tilde{\theta}_j^B\|_p \leq cC 2^{-aj}$ by the bound (59). Therefore,

$$
\begin{aligned}
R_{\text{hi}} &\leq c \sum_{j \geq j_2} 2^{jsq} 2^{-jaq} C^q + cN^{-q/2} \\
&\leq cC^q 2^{-(a-s)qj_2} + cN^{-q/2} \\
&\leq cC^{(1-2r')q} N^{-r'q} \log^{r'q} N + cN^{-q/2}.
\end{aligned}
\tag{68}
$$

To complete the proof, combine the bounds in (64), (65) and (68). □

*Discretization bias.* The risks (60), (61) and (62) all quantify errors around the discretized coefficients $\tilde{\theta}$. To control the difference in the risk norm between these coefficients and the true coefficients $\theta$, define $\Delta = (\frac{1}{q} - \frac{1}{p})_+$ as in the proof of (57). Using the bound (58) and the property $r'' = a - s - \Delta$, it follows that

$$
\begin{aligned}
\sum_{j=L-1}^{J-1} 2^{sqj} \|\tilde{\theta}_j - \theta_j\|_q^q &\leq \sum_{j=L-1}^{J-1} 2^{(s+\Delta)qj} \|\tilde{\theta}_j - \theta_j\|_p^q \\
&\leq cC^q 2^{-\bar{\alpha}qJ} \sum_{j=L-1}^{J-1} 2^{(\bar{\alpha} - r'')qj}.
\end{aligned}
\tag{69}
$$

If $\bar{\alpha} \leq r''$, the expression is bounded by $cC^q J^{\lambda'} 2^{-\bar{\alpha}qJ}$ since $\lambda' = 1$ if and only if $\bar{\alpha} = r''$. On the other hand, if $\bar{\alpha} > r''$, the sum in (69) is geometrically increasing, and so the expression is of order $cC^q 2^{-r''qJ}$. Since $r''' = \min(r'', \bar{\alpha})$, all cases are combined in the bound

$$
\sum_{j=L-1}^{J-1} 2^{sqj} \|\tilde{\theta}_j - \theta_j\|_q^q \leq cC^q N^{-r'''q} (\log N)^{\lambda'}.
\tag{70}
$$

*Completing the proof.* To complete the proof of Theorem 2, we combine the bounds (60), (61), (62) and (70). For $r' \leq r$, this gives the required result, with $\lambda = \nu$.

For $r' > r$, we have an additional term from (62), proportional to $C^{(1-2r')q} \times (N/\log N)^{-r'q}$. By elementary manipulation, we have

$$C^{(1-2r')}(N/\log N)^{-r'}$$
$$\leq \begin{cases} C^{(1-2r)}N^{-r}, & \text{if } (C^2 N)^{r'-r} \geq (\log N)^{r'}, \\ N^{-1/2}(\log N)^{(1/2)r'(1-2r)/(r'-r)}, & \text{if } (C^2 N)^{r'-r} < (\log N)^{r'}. \end{cases}$$

It follows that the bound (35) holds for this case also, setting $\lambda = \max\{\nu, \frac{1}{2}qr' \times (1-2r)/(r'-r)\}$. This completes the proof of Theorem 2.

The corresponding theorem with periodic boundary conditions on the functions and the wavelet decompositions is also true, and can be proved by a simplified version of the same approach, without any need for preconditioning or for the consideration of boundary coefficients.

To prove the result (36), we use exactly the same argument as above, but there is no need to include the discretization error or the error due to levels of the transform with $j \geq J$.

## 7. Further proofs and remarks.

7.1. *Proof of Proposition* 1. Note first that, for $q > 1$ and values of the other parameters such that $r_{q,\sigma}$ can be shown to be equivalent to the corresponding function Besov norm, the result follows from the embedding of the Besov space $B_q^{\sigma,q}$ in the Sobolev space $W_q^r$ consequent on results in Section 3.2 of [47]. However, to deal explicitly with all parameter values and with our boundary construction, we give an argument that does not use this embedding.

If $q \leq 1$, the function $\psi_{jk}^{(r)}$ has support of length at most $(2S-1)2^{-j/2}$ and maximum absolute value bounded by $c2^{j/2}2^{rj}$. Therefore, for all $j$ and $k$,

$$\int |\psi_{jk}^{(r)}|^q \leq c2^{-j/2}2^{jq/2}2^{rjq} \leq c2^{sjq}.$$

A direct calculation using the property $|x+y|^q \leq |x|^q + |y|^q$ now shows that

$$\int_0^1 |f^{(r)}(t)|^q \, dt$$

(71)
$$\leq \sum_{k \in \mathcal{K}_{L-1}} |\theta_{L-1,k}|^q \int_0^1 |\phi_{L,k}^{(r)}|^q + \sum_{j=L}^{\infty} \sum_{k=0}^{2^j-1} |\theta_{jk}|^q \int_0^1 |\psi_{jk}^{(r)}|^q$$

$$\leq c\|\theta_{L-1}\|_q^q + c\sum_{j=L}^{\infty} 2^{sjq}\|\theta_j\|_q^q$$

$$\leq cr_{q,\sigma}(f).$$

For $1 < q \leq 2$, we consider the interior and boundary cofficients of $f$ separately. Define

$$f_L = \sum_{k\in\mathcal{K}_{L-1}} \theta_{L-1,k}\phi_{Lk}, \qquad f_I = \sum_{j=L}^{\infty}\sum_{k\in\mathcal{K}_j^I} \theta_{jk}\psi_{jk}, \qquad f_B = \sum_{j=L}^{\infty}\sum_{k\in\mathcal{K}_j^B} \theta_{jk}\psi_{jk}.$$

We have, first,

$$\int_0^1 |f_L^{(r)}(t)|^q\,dt \leq \sup_t |f_L^{(r)}(t)|^q$$

(72)
$$\leq \|\theta_{L-1}\|_\infty^q \sup_t \left(\sum_{k\in\mathcal{K}_{L-1}} |\phi_{Lk}^{(r)}(t)|\right)^q$$

$$\leq c\|\theta_{L-1}\|_q^q.$$

Now consider $f_I$. Let $\chi_{jk}(x)$ be the indicator function of the interval $[2^{-j}k, 2^{-j}(k+1)]$. By Theorem 2 of Chapter 6 of [37], using the fact that $\frac{1}{2}q \leq 1$,

$$\int_0^1 |f_I^{(r)}(t)|^q\,dt \leq c\int_0^1 \left|\sum_{j=L}^{\infty}\sum_{k\in\mathcal{K}_j^I} 2^j 2^{2jr}\theta_{jk}^2 \chi_{jk}(x)\right|^{q/2} dx$$

(73)
$$\leq c\sum_{j=L}^{\infty}\sum_{k\in\mathcal{K}_j^I}\int_0^1 |2^j 2^{2jr}\theta_{jk}^2\chi_{jk}(x)|^{q/2}\,dx$$

$$= c\sum_{j=L}^{\infty} 2^{qj/2} 2^{qjr} 2^{-j}\|\theta_j\|_q^q$$

$$\leq c\sum_{j=L}^{\infty} 2^{sqj}\|\theta_j\|_q^q.$$

Finally, we consider the contribution from the boundary wavelets. Let $\mathcal{S}_j$ be the union of the supports of the boundary wavelets at level $j$. Define $\mathcal{T}_j = \mathcal{S}_j \setminus \mathcal{S}_{j+1}$ and $a_j = 2^{(r+1/2)j}\|\theta_j\|_q$. Then

$$\sum_{k\in\mathcal{K}_j^B} |\theta_{jk}\psi_{jk}^{(r)}| \leq \max_{\ell\in\mathcal{K}_j^B}|\theta_{j\ell}| \sum_{k\in\mathcal{K}_j^B} |\psi_{jk}^{(r)}(t)|$$

(74)
$$\leq c\|\theta_j\|_q 2^{(r+1/2)j} I[\mathcal{S}_j] = c a_j I[\mathcal{S}_j].$$

Using (74), the nesting properties of the $\mathcal{S}_j$, and the property that $|\mathcal{T}_\ell| = c2^{-\ell}$, we obtain, applying Hölder's inequality,

(75)
$$\int_{\mathcal{T}_\ell} |f_B^{(r)}(t)|^q \, dt \leq c \int_{\mathcal{T}_\ell} \left( \sum_{j=L}^{\infty} a_j I[\mathcal{S}_j] \right)^q = c \int_{\mathcal{T}_\ell} \left( \sum_{j=L}^{\ell} a_j \right)^q$$
$$= c2^{-\ell} \left( \sum_{j=L}^{\ell} a_j \right)^q \leq c 2^{-\ell/2} \sum_{j=L}^{\ell} 2^{-j/2} a_j^q.$$

Using the bound (75) for each term, we can now conclude that

$$\int_0^1 |f_B^{(r)}(t)|^q \, dt = \sum_{\ell=L}^{\infty} \int_{\mathcal{T}_\ell} |f_B^{(r)}(t)|^q \, dt$$
$$\leq c \sum_{\ell=L}^{\infty} \left\{ 2^{-\ell/2} \sum_{j=L}^{\ell} 2^{-j/2} a_j^q \right\}$$

(76)
$$= c \sum_{j=L}^{\infty} \left( 2^{-j/2} a_j^q \sum_{\ell=j}^{\infty} 2^{-\ell/2} \right)$$
$$= c \sum_{j=L}^{\infty} 2^{-j} a_j^q \leq c \sum_{j=L}^{\infty} 2^{sqj} \|\theta_j\|_q^q.$$

Using the bounds (72), (73) and (76), we can conclude that

$$\int_0^1 |f^{(r)}(t)|^q \, dt \leq c \left( \int |f_L^{(r)}|^q + \int |f_I^{(r)}|^q + \int |f_B^{(r)}|^q \right) \leq c r_{q,\sigma}(f),$$

completing the proof.

7.2. *How much is lost without bounded shrinkage?* The Besov space adaptivity result Theorem 1 provides a context within which the importance of the bounded shrinkage property can be assessed. Construct a coefficient array $\theta$ by setting $\theta_{3,0} = C 2^{-3a}$ and all other $\theta_{jk} = 0$. Then for any $p$, $\rho$ and $\alpha$, $\theta$ will be a member of the Besov sequence space $b_{p,\rho}^\alpha(C)$. (Obviously the choice of level 3 to have a single nonzero element is arbitrary, and any other fixed position can be used.)

Let $Z = N^{1/2} Y_{3,0}$ and $\mu = N^{1/2} \theta_{3,0}$. Use the mixture prior $(1-w)\delta_0 + w\gamma$ for $\mu$, setting $\gamma$ to be the normal $N(0, \tau^2)$ density. Whatever the value of $w$, the posterior median function then has the property

$$|\hat{\mu}(z,w)| \leq \lambda_\tau |z| \quad \text{where } \lambda_\tau = \tau^2/(1+\tau^2),$$

and the same inequality holds for the posterior mean. Since $\mu > 0$, whether $w$ is a fixed or random weight, we have

$$(77) \qquad E(\hat{\mu} - \mu)^2 \geq E(\hat{\mu} - \mu)^2 I[Z < \mu] \geq \tfrac{1}{2}(1 - \lambda_\tau)^2 \mu^2.$$

Multiplying both sides of (77) by $N^{-1}$ shows that the the mean square error risk satisfies

$$R_{N,2,0}(\theta) \geq E(\hat{\theta}_{3,0} - \theta_{3,0})^2 \geq \tfrac{1}{2}(1 - \lambda_\tau)^2 \theta_{3,0}^2 = \tfrac{1}{2}(1 - \lambda_\tau)^2 2^{-6a} C^2,$$

which does not tend to zero as $N \to \infty$. Hence, the maximum risk over any Besov sequence class does not diminish as $N$ increases, and no adaptivity result of the type given in Theorem 1 can be proved.

In the case where $\gamma$ has tails asymptotic to $\exp(-c|t|^\lambda)$ for some $\lambda \in (1, 2)$, it can be shown that, at least for large $\mu$, $|\hat{\mu} - \mu| \geq c\mu^{\lambda-1}$. Consideration of the same counterexample as above then demonstrates that, again, however the weight is chosen, the maximum risk over the Besov sequence space is bounded below by a multiple of $N^{-2+\lambda}$. For large $\alpha$ this will dominate the rate in Theorem 1, and thus the assumption that $\gamma$ has tails at least as heavy as exponential cannot essentially be relaxed without restricting or removing the adaptivity demonstrated by the theorem.

7.3. *Results for the posterior mean.* If the estimation is conducted using the posterior mean rather than a strict thresholding rule, the results of Theorems 1 and 2 still hold for $1 < q \leq 2$, since the bounds of Theorem 3 apply in this case. For smaller values of $q$ in the single sequence case, Johnstone and Silverman [33] show that the failure of the posterior mean to be a strict thresholding rule has a substantive effect on the overall risk. However, their counterexample does not unequivocally settle the question of the behavior of the posterior mean estimator in the wavelet case. The possible extension or modification of the theorems for the posterior mean estimator for $q \leq 1$ is a topic for future investigation.

**Acknowledgments.** I. M. Johnstone is very grateful for the hospitality of the University of Bristol, the "Nonparametric Semester" of Institut Henri Poincaré, Paris and the Australian National University, where parts of the work on this paper were carried out. Similarly, B. W. Silverman is grateful to the Department of Statistics at Stanford University and the Center for Advanced Study in the Behavioral Sciences at Stanford. Versions of this work were presented in B. W. Silverman's Special Invited Paper in 1999 and I. M. Johnstone's second Wald lecture in 2004. We very gratefully acknowledge Ed George for his sustained encouragement and intellectual generosity over the period of this work, and the referees for their detailed and helpful comments.

DEPARTMENT OF STATISTICS
STANFORD UNIVERSITY
STANFORD, CALIFORNIA 94305-4065
USA
E-MAIL: imj@stat.stanford.edu

ST. PETER'S COLLEGE
OXFORD OX1 2DL
UNITED KINGDOM
E-MAIL: bernard.silverman@spc.ox.ac.uk